\documentclass[]{amsart}

\usepackage[utf8]{inputenc}
\usepackage[T1]{fontenc}
\usepackage{tikz-cd}
\usepackage[margin=1.5in]{geometry}
\usepackage{graphicx}
\usepackage{adjustbox}
\usepackage{fancyhdr}
\usepackage{mathrsfs}
\usepackage{xcolor}
\usepackage{enumerate}
\usepackage{listings}
\usepackage{float}
\usepackage{hyperref}
\usepackage{pgffor}
\usepackage{cite}
\usepackage{tikz,xifthen}
\usepackage{pgf}
\usepackage{subcaption}
\usepackage{boldline}
\usepackage{mathtools}
\usepackage{parskip}
\usepackage{tabularray}
\usepackage{lscape}
\UseTblrLibrary{diagbox}
\usepackage{nicematrix}

\usepackage{verbatim}

\usetikzlibrary{arrows,positioning,decorations.markings,calc,decorations.pathmorphing}

\tikzset{parallel arrow/.style={latex-,
		shorten >=2mm, shorten <=2mm, 
		decoration={sl,raise=1cm},decorate}}

\title{Graph planar algebra embeddings and infinite depth subfactors}
\author{Dietmar Bisch}
\author{Julio C{\'a}ceres}

\address{Department of Mathematics, Vanderbilt University, 1326 Stevenson Center, Nashville, TN 37240, USA}
\email{dietmar.bisch@vanderbilt.edu}
\email{julio.e.caceres.gonzales@vanderbilt.edu}

\thanks{D.B. and J.C. were supported in part by US ARO grant W911NF-23-1-0026}

\tikzset{
	shaded/.style = {fill=red!10!blue!20!gray!30!white},
	Tbox/.style = {circle, draw, very thick, fill=white, opaque},
	PAdefn/.style = {scale=.7,baseline},
}

\newcommand{\ra}{\rightarrow}
\DeclareMathOperator{\tr}{tr}
\DeclareMathOperator{\ch}{ch}

\DeclareMathOperator{\ID}{id}

\DeclareMathOperator{\Int}{int}

\DeclareMathOperator{\Rot}{Rot}

\newcommand{\ve}{\varepsilon}
\newcommand{\II}{II\textsubscript{1}}

\newcommand{\R}{\mathbb{R}}
\newcommand{\N}{\mathbb{N}}

\newcommand{\C}{\mathbb{C}}
\newcommand{\Z}{\mathbb{Z}}

\newcommand{\T}{\mathcal{T}}
\newcommand{\V}{\mathcal{V}}

\newcommand{\E}{\mathcal{E}}
\newcommand{\sst}{\subset}
\newcommand{\bsum}{\displaystyle \sum}

\usepackage{expl3}
\ExplSyntaxOn
\cs_new_eq:NN \Repeat \prg_replicate:nn
\ExplSyntaxOff

\allowdisplaybreaks

\theoremstyle{definition}
\newtheorem{defi}{Definition}[section]
\newtheorem{rem}[defi]{Remark}
\newtheorem{nota}[defi]{Notation}
\newtheorem{exa}[defi]{Example}
\theoremstyle{plain}
\newtheorem{theo}[defi]{Theorem}
\newtheorem*{theo*}{Theorem}
\newtheorem{prop}[defi]{Proposition}
\newtheorem{lem}[defi]{Lemma}
\newtheorem{coro}[defi]{Corollary}
\numberwithin{equation}{section}

\begin{document}

\begin{abstract} 
Subfactors of the hyperfinite II$_1$ factor with ``exotic'' properties 
can be constructed from
nondegenerate commuting squares of multi-matrix algebras. We show that 
the subfactor planar algebra of these commuting square subfactors 
necessarily embeds into Jones' graph planar algebra associated to 
one of the inclusion graphs in the commuting square. This leads to 
a powerful obstruction for the standard invariant of the subfactor, 
and we use it to give
an example of a hyperfinite subfactor with Temperley-Lieb-Jones
standard invariant and index $\frac{5+ \sqrt{13}}{2}$, i.e. the index
of the Haagerup subfactor. We are led to a conjecture pertaining 
to Jones indices of irreducible, hyperfinite subfactors.
\end{abstract}

\maketitle
	


\section{Introduction}

Vaughan Jones initiated the theory of subfactors in his breakthrough
paper \cite{Jones1983} in which he introduced the notion of 
\textit{index} for an inclusion of II$_1$ factors $N \subset M$. 
He showed that the surprising result that the index is quantized 
when $\le 4$. Jones also gave an easy construction of
hyperfinite subfactors whose index can attain any real number $> 4$.
This construction was based on the fact that an amplification of
the hyperfinite II$_1$ factor $R$ is isomorphic to $R$ itself and 
results in reducible subfactors $N \subset M$ with two-dimensional
relative commutant $N' \cap M$. Jones asked in \cite{Jones1983} what 
the set of indices of \emph{irreducible} subfactors of the 
hyperfinite II$_1$ factor is, i.e. those subfactors with one-dimensional 
relative commutant. 
This question remains open to this day, but small partial progress has 
been made. For instance, it is shown in \cite{bisch1994} and later 
in \cite{stojanovic2021} that rational, non-integer numbers can be 
Jones indices of irreducible, hyperfinite subfactors. Moreover,
it is known that some of these are accumulation points of 
indices of irreducibe, hyperfinite subfactors (e.g. $4.5$
by \cite{bisch1994} and \cite{schou2013commuting}). 

Indices of finite depth subfactors of the hyperfinite II$_1$ factor
are of course algebraic integers, since they are roots of
the characteristic polynomial of the principal graphs of the subfactor.
So, in order to construct irreducible, hyperfinite subfactors with
more interesting Jones indices, one has to turn to infinite depth 
subfactors. 

A classical method to construct hyperfinite subfactors is via
(nondegenerate) commuting squares of multi-matrix algebras
(see for instance \cite{goodman2012coxeter}). A particularly
interesting class of
examples is obtained from so-called {\it spin model commuting
squares} that are built from complex Hadamard matrices. 
These are the easiest commuting squares, and they yield irreducible,
hyperfinite subfactors with integer index. There is at least one
such spin model commuting square subfactor for every integer, and
probably many more. It has been a long-standing
problem to determine the standard invariant of these subfactors, and
Jones developed {\it planar algebras} as a tool to compute the
``quantum symmetries'' (or standard invariant) encoded by these 
subfactors \cite{jones1999planar}. 
Despite partial progress (see e.g. \cite{burstein2015}, 
\cite{montgomery2023}), these subfactors remain mysterious. More
generally, computing the standard invariant of a subfactor that
is obtained from iterating the Jones basic construction of a nondegenerate
commuting square of multi-matrix algebras remains to be a formidable task.
While, in theory, the problem is solved by Ocneanu compactness 
\cite{evanskawahigashi1998}, in practice, the computation
surmounts the computational power we have. The problem is hard even in 
the case of low index subfactors, such as the
Haagerup subfactor, which is the finite depth subfactor with 
smallest index $>4$. It has index $\frac{5+\sqrt{13}}{2}$ 
\cite{asaedahaagerup1999}.

We establish in this paper a criterion that allows us to determine the 
standard invariant of a commuting square subfactor without having to do any
computations. Using the classification of standard invariants
of subfactors with index between $4$ and $5.25$ (see e.g. 
\cite{afzaly2015classification}), we
accomplish this task by proving a graph planar algebra embedding theorem 
that allows us to avoid explicit computations for the standard
invariant of commuting square subfactors. The key result in our 
strategy is theorem 1.2 in \cite{grossman2018extended} 
that establishes that the standard invariant of a finite depth subfactor
embeds into the graph planar algebra of a finite, connected and 
bipartite graph precisely when the graph is a fusion graph of a left 
$\mathcal{C}$-module C$^*$-category with 
respect to a canonical object, where $\mathcal{C}$ is the unitary fusion
category arising from the finite depth subfactor .
We show that the planar algebra of the horizontal commuting square 
subfactor must embed, as a planar algebra, into the graph planar 
algebra of the first vertical
inclusion graph appearing in the commuting square. This result is true
even if the commuting square subfactor has infinite depth. The proof
consists of showing that the action of the planar operad of shaded
planar tangles as defined by Jones for a graph planar algebra is
precisely the action of the operad defined from the subfactor when
restricted to the subspaces given by the higher relative commutants.
This can be done explicitly by using a Pimsner-Popa basis and 
the fact that all planar tangles are generated by a certain finite
family of explicit tangles \cite{kodiyalam2004jones}. Our result
is inspired by a theorem of Jones and Penneys that shows that
the subfactor planar algebra of a finite depth subfactor embeds
into the graph planar algebra of its principal graphs 
\cite{jones2011embedding}. Variations of graph planar algebra
embedding theorems were also shown in 
\cite{morrisonwalker2010}, \cite{coleshuston2021},  
\cite{chen2024stardardlattices},
\cite{montgomery2023}, and some non-embedding results obtained
via planar algebra techniques can be found in \cite{peters2010planar}.

In a forthcoming article, we will establish that every Jones index 
of an irreducible,
hyperfinite, finite depth subfactor with index $\in \left[ 4,5 \right]$ 
is also the index of a hyperfinite subfactor with Temperley-Lieb-Jones
standard invariant (and hence $A_\infty$ principal graphs). We call
such subfactors $A_\infty$ subfactors. Based on our work, we conjecture
that every Jones index of an irreducible, hyperfinite subfactor can
in fact be realized by an $A_\infty$ hyperfinite subfactor. We have shown
that this conjecture also holds true if the index is $3 + \sqrt{5}$
\cite{bischcaceres2}. Since there exists a hyperfinite
infinite depth subfactor with Fuss-Catalan standard invariant with 
this index \cite{bisch1997intermediate},
a more involved argument to show the existence of an
A$_\infty$ subfactor with index  $3 + \sqrt{5}$ is required, 
see \cite{bischcaceres2} for more details.

The classification of standard invariants of low index subfactors
shows that there are only a few Jones indices between 
$4$ and $3 + \sqrt{5}$ where finite depth subfactors occur. All 
irreducible, infinite depth (hyperfinite) subfactors with indices 
in $\left[4, 5.25\right]$ except 
these few, must be $A_\infty$ subfactors.
We summarize the Jones indices where finite depth subfactors exist
in Table \ref{indices}.

\begin{table}[ht]
    \centering
    \begin{NiceTabular}{|c|c|c|}
            \CodeBefore
            \rowcolors{1}{blue!15}{}
            \Body
            \hline
            Index & \# of subfactors & Name\\
            $\frac{1}{2}(5+\sqrt{13})$& 2& Haagerup\\
            $\approx 4.37720$& 2&Extended Haagerup\\
            $\frac{1}{2}(5+\sqrt{17})$& 2& Asaeda-Haagerup\\
            $3+\sqrt{3}$& 2& 3311\\
            $\frac{1}{2}(5+\sqrt{21})$& 2& 2221\\
            $5$& 7 & -\\
            $\approx 5.04892$& 2& $\mathfrak{su}(2)_5$ and $\mathfrak{su}(3)_4$\\
            $3+\sqrt{5}$& 11 & -\\
            \hline
        \end{NiceTabular}
    \caption{Indices of hyperfinite finite depth subfactors between $4$ and $5.25$}
    \label{indices}
\end{table}

Hyperfinite subfactors can also be constructed using planar algebra 
techniques. In \cite{peters2010planar} and 
\cite{bigelow2012extendedHaagerup}, the 
authors construct the Haagerup and extended Haagerup subfactors by 
locating their planar algebras inside the \emph{graph planar algebra} 
of their principal graphs. This method will always work for finite
depth subfactors by the Jones-Penneys result \cite{jones2011embedding}
mentioned above, at least in principle.

Let us describe in more detail the sections of this paper.
Recall that 
a \emph{commuting square} is an inclusion of four finite dimensional 
von Neumann algebras \begin{equation*}
		\begin{array}{ccc}
		A_{1,0}&\stackrel{\text{L}}{\sst} &A_{1,1}\\
		\rotatebox[origin=c]{90}{$\sst$}_K& &\rotatebox[origin=c]{90}{$\sst$}_H\\
		A_{0,0}&\stackrel{\text{G}}{\sst} &A_{0,1}
		\end{array}
\end{equation*}
as above, with a faithful trace on $A_{1,1}$ so that $A_{1,0}$ and $A_{0,1}$ are orthogonal modulo their intersection $A_{0,0}$ with respect to the inner product defined by the trace. We iterate the basic construction for the inclusions $A_{i,j}\sst A_{i,j+1}$ and $A_{i,j}\sst A_{i+1,j}$. Then, if every inclusion is connected, the commuting square is nondegenerate and the trace has the Markov property, $A_{0,\infty}\sst A_{1,\infty}$ will be an irreducible hyperfinite subfactor with index $\|K\|^2$, where $A_{n,\infty}=\left(\cup_k A_{n,k}\right)''$ ($n=0,1$) and $K$ is the inclusion graph of $A_{0,0}\sst A_{1,0}$. We call subfactors constructed in this manner \emph{commuting square subfactors}. 
We summarize the results we need regarding commuting squares 
in Section 2 (see \cite{jones1997introduction} and 
\cite{goodman2012coxeter}). 

In \cite{jones2011embedding}, the authors describe inclusions of finite 
dimensional C$^*$-algebras in terms of loop algebras, particularly those 
coming from the basic construction. Jones and Penneys present formulas for 
the Jones projections, conditional expectations and the Pimsner-Popa basis 
for these inclusions. We give a brief overview of this work in Section 3 
and a recall notions from the theory of planar algebras in Section 4.

In Section 5 we then prove our graph planar algebra embedding theorem, 
which generalizes the main theorem in \cite{jones2011embedding}:
\begin{theo}\label{maint}
    Let $P_\bullet$ be the subfactor planar algebra associated to the commuting square subfactor $A_{0,\infty}\sst A_{1,\infty}$ and let $G_\bullet$ the graph planar algebra associated to the inclusion graph of $A_{0,0}\sst A_{1,0}$. There is an embedding of planar algebras $\psi:P_\bullet\ra G_\bullet$.
\end{theo}

This theorem applies to any commuting square subfactor, regardless of 
whether it has finite or infinite depth. We use it to show the
existence of a hyperfinite subfactor with Jones index $\frac{5+\sqrt{13}}{2}$
and $A_\infty$ principal graphs.

\section{Subfactors and commuting squares}    
In this section we recall one of the main tools to construct 
(irreducible) hyperfinite subfactors with finite Jones index, namely
commuting squares. 
Commuting squares allow us to determine the index of a subfactor 
$N\sst M$ when it is approximated via finite-dimensional C*-algebras 
in the following way
	\[ \begin{array}{ccc}
	N&\sst &M\\
	\rotatebox[origin=c]{90}{$\sst$}& &\rotatebox[origin=c]{90}{$\sst$}\\
	A_{n+1}&\sst &B_{n+1}\\
	\rotatebox[origin=c]{90}{$\sst$}& &\rotatebox[origin=c]{90}{$\sst$}\\
	A_{n}&\sst &B_{n}
	\end{array} \]
	where $N=\left(\bigcup_{n\geq 0} A_n\right)''$ and 
$M=\left(\bigcup_{n\geq 0} B_n\right)''$, see \cite{pimsner1986entropy}. 
We refer to \cite{jones1997introduction} for proofs of most of the results 
in this section. 
	
	\begin{defi}
		Let $A_0\sst B_0$, $A_1\sst B_1$ be finite von Neumann algebras such that
		\begin{equation}\label{csquare}
		\begin{array}{ccc}
		A_{1}&\sst &B_{1}\\
		\rotatebox[origin=c]{90}{$\sst$}& &\rotatebox[origin=c]{90}{$\sst$}\\
		A_{0}&\sst &B_{0}
		\end{array}
		\end{equation}
		and $\tr$ is a faithful normal trace in $B_1$. Then (\ref{csquare}) is called a \emph{commuting square} if $E_{B_0}E_{A_1}=E_{A_0}$ where $E_{B_0},E_{A_1},E_{A_0}$ are the unique trace-preserving conditional expectation.
	\end{defi}

	\begin{prop}
		Consider a commuting square as in (\ref{csquare}) of finite-dimensional C*-algebras  and let $\tr$ be a Markov trace for $B_0\sst B_1$. Let $B_2=\langle B_1,e_{B_0}\rangle$ be the basic construction and let $A_2=\{A_1,e_{B_0}\}''$. Then
		\begin{equation*}
		\begin{array}{ccc}
		A_{2}&\sst &B_{2}\\
		\rotatebox[origin=c]{90}{$\sst$}& &\rotatebox[origin=c]{90}{$\sst$}\\
		A_{1}&\sst &B_{1}
		\end{array}
		\end{equation*}
		is also a commuting square.
	\end{prop}
	
	\begin{defi}
		Consider a commuting square
		\begin{equation}\label{csfindim}
		\begin{array}{ccc}
		A_{1,0}&\stackrel{\text{L}}{\sst} &A_{1,1}\\
		\rotatebox[origin=c]{90}{$\sst$}_K& &\rotatebox[origin=c]{90}{$\sst$}_H\\
		A_{0,0}&\stackrel{\text{G}}{\sst} &A_{0,1}
		\end{array}
		\end{equation}
		of finite-dimensional C*-algebras with a faithful normal trace on $A_{1,1}$ and inclusion graphs $G$, $H$, $K$ and $L$. We say that (\ref{csfindim}) is a \emph{nondegenerate} or \emph{symmetric} commuting square if $GH=KL$ and $HL^t=G^tK$.
	\end{defi}
	
	It can be shown that whenever we have a nondegenerate commuting square, the inclusion $A_0\sst A_1\sst A_2$ from the previous proposition is isomorphic to the basic construction $A_0\sst A_1\sst \langle A_1, e_{A_0}\rangle$. 
	 
	\begin{prop}
		Suppose (\ref{csfindim}) is a nondegenerate commuting square with connected graphs $G$, $H$, $K$ and $L$. Then $\|H\|=\|K\|$ and $\|G\|=\|L\|$. Moreover, $\tr$ is the Markov trace for $A_{1,0}\sst A_{1,1}$, $\tr|_{A_{0,1}}$ is the Markov trace for $A_{0,0}\sst A_{0,1}$, $\tr$ is the Markov trace for $A_{0,1}\sst A_{1,1}$ and $\tr|_{A_{1,0}}$ is the Markov trace for $A_{0,0}\sst A_{1,0}$.
	\end{prop} 
	
	\begin{proof}
		See Corollary 5.3.4 part (c) in \cite{jones1997introduction}.
	\end{proof}
	
	Thus, whenever we have a nondegenerate commuting square we can iterate the basic construction in the horizontal and vertical directions. From this, we obtain a lattice of the form
	\begin{equation}\label{infgrid}
	\begin{array}{ccccccccc}
	A_{\infty,0} & \sst & A_{\infty,1} & \sst & A_{\infty,2} & \sst & \cdots & \sst & A_{\infty,\infty} \\
	\cup &  & \cup &  & \cup &  &  &  & \cup \\
	\vdots &  & \vdots &  & \vdots &  &  &  & \vdots \\
	\cup &  & \cup &  & \cup &  &  &  & \cup \\
	A_{2,0} & \sst & A_{2,1} & \sst & A_{2,2} & \sst & \cdots & \sst & A_{2,\infty}  \\
	\cup &  & \cup &  & \cup &  &  &  &  \\
	A_{1,0}& \sst & A_{1,1} & \sst & A_{1,2} & \sst & \cdots & \sst & A_{1,\infty} \\
	\cup &  & \cup &  & \cup &  &  &  & \cup \\
	A_{0,0} & \sst & A_{0,1} & \sst & A_{0,2} & \sst & \cdots & \sst & A_{0,\infty} 
	\end{array}
	\end{equation}
	where $A_{n,\infty}=\left(\bigcup_{k} A_{n,k}\right)''$ and $A_{\infty,k}=\left(\bigcup_{n} A_{n,k}\right)''$ are hyperfinite von Neumann algebras. If we assume all inclusion matrices to be connected, then $A_{n,\infty}$ 
and $A_{\infty,k}$ are factors. 
	
	We will use this finite-dimensional approximation to compute 
invariants of the horizontal and vertical subfactors.

    \begin{defi}
        Let $A\sst (B,\tr)$ be an unital inclusion of finite von Neumann algebras with faithful normal trace $\tr$. A \emph{Pimsner-Popa} basis for $B$ over $A$ is a set $S=\{s\}\sst B$ which satisfies any of the following equivalent conditions:
        \begin{enumerate}
            \item $1=\sum_{s\in S} s e_A s^*$,
            \item $x=\sum_{s\in S} s E_A(s^* x)$ for all $x\in B$, and
            \item $x=\sum_{s\in S} E_A(x s)s^*$ for all $x\in B$,
        \end{enumerate}
        where $E_A:B\ra A$ is the conditional expectation determined by $\tr$ and $e_A$ is the Jones projection associated to this inclusion (see \cite{jones1997introduction}). 
    \end{defi}
    
	\begin{lem}
		Suppose we have a nondegenerate commuting square as in (\ref{csfindim}) with respect to the normalized Markov trace. Then there exists a finite set $I$, and $\{s_i:\;i\in I\}\sst A_{1,0}$, $\{f_i:\; i\in I\}\sst A_{0,0}$ such that 
		\begin{enumerate}[(a)]
			\item each $f_i$ is a projection,
			\item $E_{A_{0,0}}(s_i^* s_j)=\delta_{ij}f_i$,
			\item $\sum_{i\in I}\tr(f_i)=\|K\|^2$, and 
			\item $x=\sum_{i\in I} E_{A_{0,1}}(xs_i)s_i^*$, for all $x\in A_{1,1}$.
		\end{enumerate}
	\end{lem}
	
	\begin{proof}
		See Lemma 5.7.3 in \cite{jones1997introduction}.
	\end{proof}
	
	\begin{rem}
		Note that $\{s_i:\;i\in I\}\sst A_{1,0}$ is a Pimsner-Popa basis for $A_{1,1}$ over $A_{0,1}$ and for $A_{1,0}$ over $A_{0,0}$.
	\end{rem}
	
	\begin{coro}\label{ppbasesCS}
		Suppose $\{s_i:\; i\in I\}$ is as in the previous Lemma, then
		\begin{enumerate}[(a)]
			\item $\{s_i:\; i\in I\}$ is a Pimsner-Popa basis for $A_{0,\infty}\sst A_{1,\infty}$.
			\item $[A_{1,\infty}:A_{0,\infty}]=\|K\|^2$.
		\end{enumerate}
	\end{coro}

	\begin{proof}
		See Corollary 5.7.4 in \cite{jones1997introduction}.
	\end{proof}

    Connected unital inclusions of finite dimensional C*-algebras and finite index inclusions of \II factors are both examples of \emph{strongly Markov} inclusions (see \cite{jones2011embedding}). In this case, we have the following formula for the conditional expectation on the relative commutants.
    
    \begin{prop}\label{ppbasesrelcom}
		If $A\sst B$ is strongly Markov and $[B:A]<\infty$, the conditional expectation $E_{B'}:A'\cap B(L^2(B_n,\tr))\ra B'\cap B(L^2(B_n,\tr))$ is given by 
		\[ E_{B'}(x)=\frac{1}{[B:A]}\sum_{s\in S} sxs^* \]
		where $S$ is a Pimsner-Popa basis for $B$ over $A$. In particular, the map is independent of choice of a basis.
	\end{prop}

    \begin{proof}
		See Proposition 2.23 in \cite{jones2011embedding} or Proposition 2.7 in \cite{bisch1997bimodules}.
	\end{proof}
    
    \begin{rem}
        When $A\sst B$ is a connected unital inclusions of finite dimensional C*-algebras then $[B:A]=\|\Gamma\|^2$ where $\Gamma$ is the inclusion matrix for $A\sst B$.
    \end{rem}
    
	The next result will allow us, in principle, to compute the standard invariant from our finite-dimensional commuting squares. 
	
	\begin{theo}[Ocneanu compacntess]\label{ocneanu}
		Given a lattice of inclusions arising from a nondegenerate commuting square as in (\ref{infgrid}), we have
		\[ A'_{0,\infty}\cap A_{n,\infty}= A'_{0,1}\cap A_{n,0},\quad A'_{\infty,0}\cap A_{\infty,n}= A'_{1,0}\cap A_{0,n},\text{ for all }n\geq 0.\]
	\end{theo}
	\begin{proof}
		See Theorem 5.7.6 in \cite{jones1997introduction}.
	\end{proof}

\section{Finite-dimensional inclusions and loops}
	In this section we will use loop algebras to describe the iterated basic construction of a finite dimensional inclusion of C*-algebras. Moreover, we will be able to describe the commutants as loop algebras on the Bratteli diagram for the inclusion. Many results in this section can be found in Section 3 
of \cite{jones2011embedding}.
	
	\subsection{Loop algebras}
	We begin by introducing the notation used in \cite{jones2011embedding}. 
	\begin{nota}\label{loopnot} Let $\Gamma$ be a finite, connected, bipartite multigraph. 
		\begin{itemize}
			\item $\V_+$ denotes the set of even vertices.
			\item $\V_{-}$ denotes the set of odd vertices.
			\item $\E$ denotes the set of edges and its elements are denoted by $\varepsilon$, $\xi$ and $\zeta$.
			\item All edges are directed from even to odd vertices indicated by the source and target functions $s:\E \ra \V_+$ and $t:\E\ra \V_-$.
			\item The set of all edges traversed from odd to even vertices are $\E^*=\{\varepsilon^*,\; \varepsilon\in\E\}$ with source and target functions $s:\E^*\ra \V_-$ and $t:\E^*\ra \V_+$ such that $s(\ve^*)=t(\ve)$ and $t(\ve^*)=s(\ve)$.
			\item The dimension vector for the even vertices is a function $m_+:\V_+\ra \N$.
			\item The dimension vector for odd vertices $m_-:\V_-\ra \N$ is defined by
			\[ m_-(v)=\sum_{t(\varepsilon)=v}m_+(s(\varepsilon)) \]
		\end{itemize}
	\end{nota}
	

	\begin{defi}
		Let $G_{0,\pm}$ be the complex vector space with basis $\V_\pm$ respectively. For $n\in \N$, $G_{n,\pm}$ will denote the complex vector space with basis loops of length $2n$ on $\Gamma$ based at a vertex in $\V_\pm$ respectively.
	\end{defi}
	
	We now describe the algebraic structure of $G_{n,+}$. In a similar manner, we obtain the algebraic structure of $G_{n,-}$.
	
	\begin{nota}
		Loops in $G_{n,+}$ are denoted by $[\ve_1 \ve_2^*\cdots \ve_{2n-1}\ve_{2n}^*]$. Any time we write this, it is implied that 
		\begin{itemize}
			\item $t(\ve_i)=s(\ve_{i+1}^*)=t(\ve_{i+1})$ for all odd $i<2n$,
			\item $t(\ve_i^*)=s(\ve_i)=s(\ve_{i+1})$ for all even $i<2n$, and
			\item $t(\ve_{2n}^*)=s(\ve_{2n})=s(\ve_1)$.
		\end{itemize}
		For a loop $x=[\ve_1 \ve_2^*\cdots \ve_{2n-1}\ve_{2n}^*]\in G_{n,+}$, we define the truncated $x_{[j,k]}$ to be the $j^\text{th}$ to $k^\text{th}$ entries of $x$ for $1\leq j\leq k \leq 2n$. For example if $j$ and $k$ are odd then $x_{[j,k]} = [\ve_j \ve_{j+1}^*\cdots \ve_{k-1}^*\ve_{k}]$. In particular $x_{[1,2n]}=x$.
	\end{nota}

	\begin{defi}
		We define an antilinear map $*$ on $G_{n,+}$ by the antilinear extension of the map 
		\[ [\ve_1 \ve_2^*\cdots \ve_{2n-1}\ve_{2n}^*]^*= [\ve_{2n} \ve_{2n-1}^*\cdots \ve_{2}\ve_{1}^*]. \]
		Note that this is just traversing the loop in the opposite direction. It is clear that we can extend this notion for truncated paths $x_{[j,k]}$. We can also define a multiplication on $G_{n,+}$ by
		\[ x\cdot y = \delta_{x^*_{[n+1,2n]},y_{[1,n]}}[x_{[1,n]}y_{[n+1,2n]}]. \]
		This is $0$ if the second half of $x$ and the first half of $y$ do not coincide, otherwise, we just concatenate the first half of $x$ and the second half of $y$. Clearly, $*$ is an involution for $G_{n,+}$ under this multiplication.
	\end{defi}
	
	\begin{defi}
		Let $\widetilde{\Gamma}$ be the augmentation of the bipartite graph $\Gamma$ by adding a distinguished vertex $\star$ which is connected to each $v\in \V_+$ by $m_+(v)$ distinct edges, all oriented so they begin at $\star$. We will denote these added edges by $\eta$'s, $\kappa$'s or $\rho$'s.
	\end{defi}

	\begin{exa}
		Suppose the dimension vector for $\V_+$ only takes the value $1$. Then augmenting a graph $\Gamma$ looks as follows
		$$\Gamma=
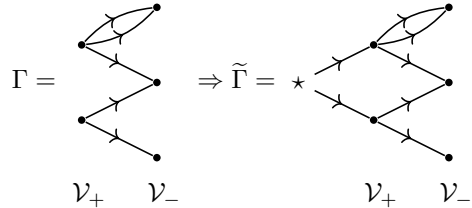
\begin{tikzpicture}[semithick,decoration={markings,mark=at position 0.5 with {\arrow{>}}},baseline=-0.5ex]
			\node[circle,fill=black,inner sep=0pt,minimum size=1mm] at (0,0) (A12) {};
			\node[circle,fill=black,inner sep=0pt,minimum size=1mm] at (-1,0.5) (A01)  {};
			\node[circle,fill=black,inner sep=0pt,minimum size=1mm] at (-1,-0.5) (A02) {};
			\node[circle,fill=black,inner sep=0pt,minimum size=1mm] at (0,1) (A11) {};
			\node[circle,fill=black,inner sep=0pt,minimum size=1mm] at (0,-1) (A13) {};
			\draw[postaction={decorate}] (A01) to [bend left =20] (A11);
			\draw[postaction={decorate}] (A01) to [bend right =20] (A11);
			\draw[postaction={decorate}] (A01) --  (A12);
			\draw[postaction={decorate}] (A02) --  (A12);
			\draw[postaction={decorate}] (A02) --  (A13);
			\node at (-0.9,-1.5) {$\V_+$};
			\node at (0.1,-1.5) {$\V_-$};
			\end{tikzpicture}\Rightarrow \widetilde{\Gamma}=
			\begin{tikzpicture}[semithick,decoration={markings,mark=at position 0.5 with {\arrow{>}}},baseline=-0.5ex]
			\node[circle,fill=black,inner sep=0pt,minimum size=1mm] at (0,0) (A12) {};
			\node[circle,fill=black,inner sep=0pt,minimum size=1mm] at (-1,0.5) (A01)  {};
			\node[circle,fill=black,inner sep=0pt,minimum size=1mm] at (-1,-0.5) (A02) {};
			\node[circle,fill=black,inner sep=0pt,minimum size=1mm] at (0,1) (A11) {};
			\node[circle,fill=black,inner sep=0pt,minimum size=1mm] at (0,-1) (A13) {};
			\node at (-2,0) (star) {$\star$};
			
			\draw[postaction={decorate}] (A01) to [bend left =20] (A11);
			\draw[postaction={decorate}] (A01) to [bend right =20] (A11);
			\draw[postaction={decorate}] (A01) --  (A12);
			\draw[postaction={decorate}] (A02) --  (A12);
			\draw[postaction={decorate}] (A02) --  (A13);
			\draw[postaction={decorate}] (star) --  (A01);
			\draw[postaction={decorate}] (star) --  (A02);
			
			\node at (-0.9,-1.5) {$\V_+$};
			\node at (0.1,-1.5) {$\V_-$};
			\end{tikzpicture}
        $$
	\end{exa} 
	
	\begin{defi}\label{towerB}
		For $n\in \Z_{\geq 0}$, let $B_n$ be the complex vector space with basis given by loops of length $2n+2$ in $\widetilde{\Gamma}$ that are of the form
		\[ [\eta_1\ve_1\ve_2^*\cdots \ve_{2n-1}\ve_{2n}^*\eta_2^*] \]
		this means that the loop starts and ends at $\star$ and remains in $\Gamma$ otherwise. We define the involution and multiplication for $B_{n}$ in a similar manner to our definition for $G_{n,+}$.
	\end{defi}
	
	\begin{rem}
		Sometimes it is convenient to write an element $[\eta_1\ve_1\ve_2^*\cdots \ve_{2n-1}\ve_{2n}^*\eta_2^*]$ in $B_n$ as $[\ell_1(\ell_2)^*]$ where $\ell_i$ are paths of length $n+1$ that start at $\star$ but remain in $\Gamma$ otherwise. With this notation, $t(\ell_1)=t(\ell_2)$ is $t(\ve_n)$ when $n$ is odd and $s(\ve_n)$ when $n$ is even. 
		
	\end{rem}

	\begin{defi}
		Consider the map $\iota:B_n\ra B_{n+1}$ defined by the linear extension of 
		\[ [\ell_1(\ell_2)^*]\mapsto \begin{cases}
		\bsum_{s(\ve)=s(\ve_n)} [\ell_1\ve\ve^*(\ell_2)^*] &n\text{ even}\\
		\bsum_{s(\ve)=t(\ve_n)} [\ell_1\ve^*\ve(\ell_2)^*] &n\text{ odd}
		\end{cases} \]
	\end{defi}

	It is not hard to see that $\iota$ is an injective $*$-homomorphism and therefore we can identify $B_n$ with its image $\iota(B_n)\sst B_{n+1}$. With this identification, we obtain a tower of algebras
	\[ B_0\sst B_1\sst B_2\sst \cdots \sst B_n\sst \cdots \] 
	Moreover, this allows us to define the product of elements in $B_m$ and $B_n$ just by including both in $B_{\max\{m,n\}}$. More precisely, if $x\in B_m$ and $y\in B_n$ with $m\leq n$  then
	\[ x\cdot y = \delta_{x^*_{[m+2,2m+2]},y_{[1,m+1]}}[x_{[1,m+1]}y_{[m+2,2n+2]}] \]
	
	\subsection{Isomorphism of towers}
	
	Let $A_0\sst A_1$ be an unital inclusion of finite-dimensional C*-algebras with Bratteli diagram $\Gamma$ and inclusion matrix $\Lambda$. Thus, if $A_0\simeq \bigoplus_{s=1}^{k} M_{\nu_i}(\C)$ and $A_1\simeq \bigoplus_{r=1}^{l} M_{\mu_i}(\C)$ then $\Lambda$ is a $k\times l$ matrix where $\Lambda_{sr}$ is the number of times the $s$-th summand of $A_0$ sits in the $r$-th summand of $A_1$. Recall that the Bratteli diagram is a bipartite graph with vertices $\V_+=\{v^+_1,\dots,v^+_k\}$ which correspond to the simple summands of $A_0$ and vertices $\V_-=\{v^-_1,\dots,v^-_l\}$ which correspond to the simple summands of $A_1$. By definition of the Bratteli diagram, we have $\Lambda_{ij}$ edges going from $v^+_i$ to $v^-_j$. Since our inclusion is unital, the dimension vectors $\vec{\nu}$ and $\vec{\mu}$ associated to $A_0$ and $A_1$, respectively, satisfy the equation $\vec{\mu}=\Lambda^t \vec{\nu} $. Suppose now that $\tr$ is a $d$-Markov trace on $A_1$. Hence we can iterate the basic construction to obtain a tower of algebras 
	\[ A_0\sst A_1\sst A_2\sst \cdots \sst A_n \sst \cdots \]
	It is a well-known fact (see \cite[Section 2.4]{goodman2012coxeter} or \cite[Section 3.2]{jones1997introduction}) that if we do the basic construction to obtain $A_1\sst A_2$, then $A_2$ will also be a finite-dimensional C*-algebra and moreover it has inclusion matrix $\Lambda^{t}$. Hence its Bratteli diagram is just a reflection of $\Gamma$. Consequently, if we iterate the basic construction, the inclusion matrix for $A_m\sst A_{m+1}$ is $\Lambda$ if $m$ is even and $\Lambda^t$ if $m$ is odd. Moreover, the dimension vector for $A_n$, denoted $\vec{\nu}^n$, is given inductively by $\Lambda \vec{\nu}^{n-1}$ if $n$ is even or $\Lambda^t \vec{\nu}^{n-1} $ if $n$ is odd. Thus
	\begin{align*}
	\nu^n_s &= \begin{cases}
	\bsum_{r=1}^{k} \nu^{n-1}_r \Lambda_{sr} &,\; $n$\text{ even}\\
	\bsum_{r=1}^{l} \nu^{n-1}_r \Lambda_{rs} &,\; $n$\text{ odd}
	\end{cases}
	\end{align*}

	Note that we can think of $\vec{\nu}^n$ as a function from $\V_+$ to $\N$ if $n$ is even and from $\V_{-}$ to $\N$ if $n$ is odd. Hence $\vec{\nu}^0$ will play the same role as $m_+$ in \ref{loopnot}, which allows us to define $\widetilde{\Gamma}$. From $\widetilde{\Gamma}$ we obtain the tower of loop algebras $(B_n)_{n\geq 0}$, this tower is isomorphic to the tower of algebras $(A_n)_{n\geq 0}$ obtained from the basic construction (see \cite[Section 5.4]{jones1997introduction} or \cite[Section 3.2]{jones2011embedding}).
 	
    Recall that the tower $(A_n)_{n\geq 0}$ has additional structure given by the Markov traces, conditional expectations and Jones projections. The next step is to determine how these are expressed in the loop algebra representations $(B_n)_{n\geq 0}$. 
	
	\begin{defi}\label{markovTracev}
		Let $\lambda^i=(\lambda^i_s)$ be the Markov trace (column) vector for $A_i$ for $i=0,1$ such that
		\[ \sum_s \nu_s^0\lambda_s^0 = 1 =\sum_s \nu_s^1\lambda_s^1\]	
		$\lambda^i_s$ is the trace of a minimal projection in the $s$-th simple summand of $A_i$ for $i=0,1$. Since the trace on $A_1$ restricts to the trace on $A_0$ we must have $\Lambda \lambda^1 = \lambda^0$. 
	\end{defi}
	
	\begin{rem}
		Since $\tr$ is a $d$-Markov trace, we also have $\Lambda\Lambda^t\lambda^0=d^{-2}\lambda^0$ and $\Lambda^t\Lambda \lambda^1=d^{-2}\lambda^1$, and $d=\sqrt{\|\Lambda^t\Lambda\|}=\sqrt{\|\Lambda\Lambda^t\|}$. This implies that the trace vector $\lambda^n$ for $A_n$ is given by $\lambda^n=d^{-2}\lambda^{n-2}$ for all $n\geq 2$. In particular we have
		\[ \lambda^n = \begin{cases}
		d^{-n} \lambda^0&,\; n\text{ even}\\
		d^{-n+1} \lambda^1&,\; n\text{ odd}
		\end{cases} \]
	\end{rem}

	Since $\V_+=\{v^+_1,\dots,v^+_k\}$ and $\V_-=\{v^-_1,\dots,v^-_l\}$, if $n$ is even, we can think of $\lambda^n$ as function from $\V_+$ to $\R_+$ by setting $\lambda^n(v^+_i):=\lambda^n_i$. Similarly, if $n$ is odd, we can think of $\lambda^n$ as function from $\V_-$ to $\R_+$.
	
	\begin{defi}[Traces] We define a normalized trace on $B_n$ by
		\[ \tr_n([\ell_1(\ell_2)^*])=\delta_{\ell_1,\ell_2}\lambda^n(t(\ell_1)). \]
	\end{defi}

	\begin{prop}(Conditional Expectations)\label{graphCExp}
		If $x=[\eta_1\ve_1\ve_2^*\cdots \ve_{2n-1}\ve_{2n}^*\eta_2^*]\in B_n$, the $\tr_n$-preserving conditional expectation $B_n\ra B_{n-1}$ is given by 
		\[  E_{B_{n-1}}(x)=\begin{cases}
		d^{-2}\delta_{\ve_n,\ve_{n+1}}\left(\frac{\lambda^0(s(\ve_n))}{\lambda^1(t(\ve_n))}\right)[\eta_1\ve_1\ve_2^*\cdots\ve_{n-1}\ve_{n+2}^*\cdots\ve_{2n-1}\ve_{2n}^*\eta_2^*]&,\; n\text{ even}\\
		\delta_{\ve_n,\ve_{n+1}}\left(\frac{\lambda^1(t(\ve_n))}{\lambda^0(s(\ve_n))}\right)[\eta_1\ve_1\ve_2^*\cdots\ve_{n-1}^*\ve_{n+2}\cdots\ve_{2n-1}\ve_{2n}^*\eta_2^*]&,\; n\text{ odd.}
		\end{cases}\]
	\end{prop}

	\begin{rem}
		If $n$ is even and we write $x=[\ell_1 \ve_n^*\ve_{n+1}(\ell_2)^*]\in B_n$ (this means $t(\ell_1)=t(\ve_n)$ and $t(\ell_2)=t(\ve_{n+1})$) where $\ell_i$ are paths of length $n$ starting at $\star$ and remaining in $\Gamma$ otherwise, we have
		\[ E_{B_{n-1}}(x)=d^{-2}\delta_{\ve_n,\ve_{n+1}}\left(\frac{\lambda^0(s(\ve_n))}{\lambda^1(t(\ve_n))}\right)[\ell_1(\ell_2)^*]. \]
		We obtain a similar expression in the case of $n$ odd.
	\end{rem}
	
	\begin{defi}(Jones Projections)\label{JonesprojB}
		For $n\geq 1$, consider the following elements in $B_{n+1}$: If $n$ is odd, define
		\[ F_n=\bsum_{\ell}\bsum_{\substack{s(\ve_n)=t(\ell)\\s(\ve_{n+1})=t(\ell)}} \frac{d[\lambda^1(t(\ve_{n}))\lambda^1(t(\ve_{{n+1}}))]^{1/2}}{\lambda^0(s(\ve_{n}))}[ \ell\ve_{{n}}\ve^*_{{n}}\ve_{{n+1}}\ve_{{n+1}}^*\ell^*]\]
		and if $n$ is even, define
		\[ F_n=\bsum_{\ell}\bsum_{\substack{t(\ve_n)=t(\ell)\\ t(\ve_{n+1})=t(\ell)}}\frac{[\lambda^0(s(\ve_{n}))\lambda^0(s(\ve_{{n+1}}))]^{1/2}}{d\lambda^1(t(\ve_{n}))}[ \ell\ve^*_{{n}}\ve_{{n}}\ve^*_{{n+1}}\ve_{{n+1}}\ell^*] \]
		where the sum, in both cases, is taken over all paths $\ell$ of length $n$ starting at $\star$ but remaining in $\Gamma$ otherwise.
	\end{defi}

	\begin{prop}[Basic Construction]\label{bcAandB}\label{towerisom}
		For $n\in \N$, the inclusion $$B_{n-1}\sst B_n\sst (B_{n+1},\tr_{n+1},d^{-1}F_n)$$ is standard. Hence for all $k>0$, there are isomorphisms $\varphi_k:A_k\ra B_k$ preserving the trace such that $\varphi_{k+1}|_{A_k}=\varphi_k$ and $\varphi_m(f_n)=d^{-1}F_n$ for all $m>n$, where $f_n$ is the Jones projection for the inclusion $A_{n-1}\sst A_n$.
	\end{prop}

    \begin{proof}
        See Proposition 3.22 in \cite{jones2011embedding}.
    \end{proof}

	\subsection{Relative commutants of loop algebras}\label{relcomloops}
	
	We will now provide isomorphisms between the relative commutants of the tower $(B_n)_{n\geq 0}$ and the algebras $G_{n,\pm}$. To simplify notation we will represent elements of $G_{n,\pm}$ as $[\ell^\pm_1(\ell^\pm_2)^*]$ where $\ell^\pm_i$ is a path of length $n$ in $\Gamma$ starting from $\V_\pm$.
	
	\begin{prop}
		A basis for $B'_0\cap B_n$ is given by
		\[ S_{0,n}=\left\{\bsum_{t(\eta)=s(\ell_1^+)}[\eta\ell_1^+(\ell_2^+)^*\eta^*]\in B_n;\; [\ell^+_1(\ell^+_2)^*]\in G_{n,+} \right\}. \]
		A basis for $B'_1\cap B_{n+1}$ is given by 
		\[ S_{1,n+1}=\left\{\bsum_{\substack{t(\eta)=s(\ve)\\t(\ve)=s(\ell_1^-)}}[\eta\ve\ell_1^-(\ell_2^-)^*\ve^*\eta^*]\in B_n;\; [\ell^-_1(\ell^-_2)^*]\in G_{n,-} \right\}. \]
	\end{prop}
	
	\begin{defi}
		For $n\in \Z_{\geq 0}$, let $H_{n,+}=B_0'\cap B_n$, $H_{n,-}=B_1'\cap B_{n+1}$, $Q_{n,+}=A'_0\cap A_n$, and $Q_{n,-}=A'_1\cap A_{n+1}$.
	\end{defi}
	
	\begin{coro}\label{loopIsom}
		There are canonical algebra $*$-isomorphisms $\phi_{n,\pm}:H_{n,\pm}\ra G_{n,\pm}$. If $n=0$ the isomorphisms are given by
		\[ \phi_{0,+}\left(\bsum_{t(\eta)=v_+}[\eta\eta^*]\right)=[v_+]\text{  and  } \phi_{0,-}\left(\bsum_{\substack{t(\eta)=s(\ve)\\ t(\ve)=v_-}}[\eta\ve\ve^*\eta^*]\right)=[v_-].\]
		For $n>0$, the isomorphisms are given by
		\[ \phi_{n,+}\left(\bsum_{t(\eta)=s(\ell_1^+)}[\eta\ell_1^+(\ell_2^+)^*\eta^*]\right)=[\ell_1^+(\ell_2^+)^*]\text{  and  } \phi_{n,-}\left(\bsum_{\substack{t(\eta)=s(\ve)\\ t(\ve)=s(\ell_1^-)}}[\eta\ve\ell_1^-(\ell_2^-)^*\ve^*\eta^*]\right)=[\ell_1^-(\ell_2^-)^*]. \]
	\end{coro}

	\begin{rem}\label{towerisomor}
		For $n\geq 0$, $\psi_{n,+}:=\phi_{n,+}\circ \left(\varphi_n|_{H_{n,+}}\right)$ and $\psi_{n,-}:=\phi_{n,-}\circ \left(\varphi_{n+1}|_{H_{n,-}}\right)$ are isomorphisms from $Q_{n,\pm}$ to $G_{n,\pm}$. 
	\end{rem}
	
	We describe next the conditional expectation on the relative commutants using an explicit Pimsner-Popa basis for $B_1$ over $B_0$.
	
	\begin{prop}[Pimsner-Popa Basis]\label{PPbases} For each $v_+\in \V_+$, pick a distinguished $\eta_{v_+}$ with $t(\eta_{v_+})=v_+$. Set
		\begin{align*}
		S_1 &=\left\{ \left(\frac{\lambda^0(s(\ve_2))}{\lambda^1(t(\ve_2))}\right)^{1/2}\bsum_{t(\eta)=s(\ve_1)}[\eta\ve_1\ve_2^*\eta^*] :\; [\ve_1\ve_2^*]\in G_{1,+} \right\},\\
		S_2 &=\left\{ \left(\frac{\lambda^0(s(\ve_2))}{\lambda^1(t(\ve_2))}\right)^{1/2}[\eta_{s(\ve_1)}\ve_1\ve_2^*\eta_{s(\ve_2)}^*] :\; s(\ve_1)\neq s(\ve_2) \right\}.
		\end{align*}
		Then $S=S_1 \coprod S_2$ is a Pimsner-Popa basis for $B_1$ over $B_0$.
	\end{prop}
	
	\begin{coro}[Commutant Conditional Expectations]\label{loopcomCond}
		If
		\[ x = \bsum_{t(\rho)=s(\xi_1)}[\rho\xi_1\ell_1^-(\ell_2^-)^*\xi_{2n}^*\rho^*]\in B'_0\cap B_n\]
		the $\tr_n$-preserving conditional expectation $B'_0\cap B_n\ra B'_1\cap B_n$ is given by 
		$$E_{B'_1}(x)=d^{-2}\delta_{\xi_1,\xi_{2n}}\left(\frac{\lambda^0(s(\xi_1))}{\lambda^1(t(\xi_1))}\right)\bsum_{\substack{t(\rho)=s(\ve)\\ t(\ve)=t(\xi_1)}}[\rho\ve\ell_1^-(\ell_2^-)^*\ve^*\rho^*]$$
	\end{coro}

\section{Planar algebras}
	Planar algebras were introduced by Jones in \cite{jones1999planar} to describe the algebraic structure of the standard invariant of subfactors. In this section, we will give a broad description of a planar algebra. Our main interest lies in the definition of a subfactor planar algebra and a graph planar algebra. 
	
	\begin{defi}
		A \emph{planar $n$-tangle} $T$ consists of the unit disk $D$ ($=D_0$) in $\R^2$ together with a finite (possibly empty) set of disjoint subdisks $D_1$, $D_2$, $\dots$ , $D_k$ in the interior of $D$. Each disk $D_i$, $i\geq 0$, will have an even number $2n_i\geq 0$ of marked points on its boundary (with $n=n_0$). Inside $D$ there is also a finite set of disjoint smoothly embedded curves called \emph{strings} which are either closed curves or whose boundaries are marked points of the $D_i$'s. Each marked point is the boundary of some string, which meets the boundary of the corresponding disk transversally. The strings all lie in the complement of the interiors $\Int(D_i)$ of the $D_i$, $i>0$. The connected components of the complements of the strings in $\Int(D)\setminus \bigcup_{i=1}^{k} D_i$ are called \emph{regions} and are shaded black and white so that the regions whose closures meet have different shadings. The shading is part of the data of the tangle, as is the choice, at every $D_i$, $i\geq 0$, of a region whose closure meets that disk which is marked with a $\$$.
	\end{defi}
		
	Here is an example of a planar $3$-tangle:
	\[ T= 
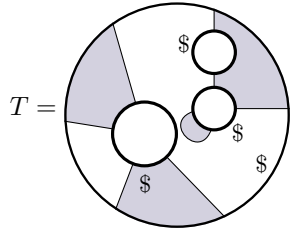
\begin{tikzpicture}[scale=.5,baseline=0]
	\clip (0,0) circle (3cm);
	
	\begin{scope}[shift=(10:1cm)] 
	\draw[shaded] (0,0)--(0:6cm)--(90:6cm)--(0,0);  
	\draw[shaded] (0,0) .. controls ++(180:2cm) and ++(-90:2cm) .. (0,0);
	\end{scope}
	
	\draw[shaded] (-150:1cm) -- (120:4cm) -- (180:4cm) -- (-150:1cm);
	\draw[shaded] (-150:1cm) -- (-120:4cm) -- (-60:4cm) -- (-150:1cm);
	
	\begin{scope}[shift=(10:1cm)] 
	\node at (0,0) [Tbox, inner sep=2mm] {};
	\node at (90:1.5cm) [Tbox, inner sep=2mm] {};
	\node at (-45:.9cm) {\small$\$$};
	\node at (115:1.9cm) {\small$\$$};
	\end{scope}
	\node at (-150:1cm) [Tbox, inner sep=3mm] {};
	\node at (-115:2cm) {\small$\$$};
	\node at (-30:2.6cm) {\small$\$$};
	
	\draw[very thick] (0,0) circle (2.98cm);
	\end{tikzpicture}	\]

	\begin{nota}
		For every disk $D_i$ in a tangle we let $\partial(D_i)=(n_i,+)$ if the region corresponding to the marked region is white and $\partial(D_i)=(n_i,-)$ if the shading is black. In particular, if we set $\partial(T)=\partial(D_0)$ then for the example above we have $\partial(T)=(3,+)$. We will say that a $k$-tangle $T$ is \emph{positive} if $\partial(T)=(n,+)$ and \emph{negative} if $\partial(T)=(n,-)$.
	\end{nota}
	
	\begin{rem}
		We are only interested in the isotopy class of a planar tangle. Later on we will see that it is useful to deform all the circles in our diagrams into rectangles to get a ``standard form'' for a tangle.
	\end{rem}
	
	\begin{defi}[Composition of tangles]
		Consider $T$ a $n$-tangle, $S$ a $n'$-tangle and $D_i$ a disk in $T$ such that $n_i=n'$. If the shading associated to the marked interval of $D_i$ coincides with the shading of the marked interval of the boundary of $S$ then we define $T\circ_i S$ to be the $n$-tangle obtained by ``gluing'' $S$ into $D_i$ such that the marked interval and shading of $D_i$ and $S$ coincide. 
	\end{defi}

	\begin{defi}[Planar algebra]
		A \emph{planar algebra} is a family $\{P_{n,\pm}\}_{n\in \N}$ of complex vector spaces together with linear maps
		\[ Z_T: \bigotimes_{i=1}^k P_{\partial(D_i)}\ra P_{\partial(T)} \]
		for every $n$-tangle $T$. 
	\end{defi}

	We require the maps $Z_T$ to only depend on the isotopy class of $T$, be independent of the ordering of the internal disks of $T$ and to be ``compatible with composition of tangles'' in the following manner. If $T$ is a $n$-tangle with $k$ internal disks $\{D^T_i\}_{i=1}^{k}$ and $S$ is a $n'$-tangle with $r$ internal disks $\{D^S_i\}_{i=1}^{r}$, then $T\circ_i S$ is an $n$-tangle with $k+r-1$ internal disks given by 
	\[ D_j= \begin{cases}
	D^T_j & 1\leq j<i\\
	D^S_{j-i+1}& i\leq j\leq i+r-1\\
	D^T_{j-r+1}& i+r\leq j\leq r+k-1
	\end{cases} \]
	We say $Z$ is compatible with the composition of tangles if the following diagram commutes:
	
	\begin{tikzcd}
	\left(\otimes_{j=1}^{i-1} P_{\partial(D^T_j)}\right)\otimes \left(\otimes_{j=1}^{r} P_{\partial(D^S_j)}\right)\otimes \left(\otimes_{j=i+1}^{k} P_{\partial(D^T_j)} \right) \arrow[rd,"Z_{T\circ_i S}"] \arrow[dd, "\ID\otimes Z_S\otimes\ID"]& \\
	& P_{\partial(T)}\\
	\otimes_{j=1}^{k} P_{\partial(D^T_j)}\arrow[ur, "Z_T"] & 
	\end{tikzcd}
	
	When the tangle $T$ has no internal disks we follow the convention that an empty tensor product is the underlying field, in this case $\C$. This means that each subspace $P_{n,\pm}$ has a distinguished subset $\{Z_T(1): T\text{  a $n$-tangle without internal disks}\}$. We call these \emph{Temperley-Lieb tangles}.
	
	We are interested in the case where $P_{n,\pm}$ is finite-dimensional for all $n\geq 0$. If $P_{0,\pm}$ is one-dimensional, there is a unique way to identify $P_{0,+}$ with $\C$ such that $Z_{\begin{tikzpicture}[scale=0.1]
		\draw[thick] (0,0) circle (1cm);
		\end{tikzpicture}}(1)=1$ and $P_{0,-}$ with $\C$ such that $Z_{\begin{tikzpicture}[scale=0.1]
		\draw[thick,shaded] (0,0) circle (1cm);
		\end{tikzpicture}}(1)=1$, where both tangles have no internal disks nor strings but have different shading. There are two scalars associated to such a planar algebra, 
	\[  \delta_1=Z_{\begin{tikzpicture}[scale=0.2]
		\draw[thick] (0,0) circle (1cm);
		\draw[shaded] (0,0) circle (0.5cm);
		\end{tikzpicture}}(1),\quad \delta_2=Z_{\begin{tikzpicture}[scale=0.2]
		\draw[thick,shaded] (0,0) circle (1cm);
		\draw[fill=white] (0,0) circle (0.5cm);
		\end{tikzpicture}}(1),\]
	where the inner circles are strings, not disks. We call $\delta_1$ and $\delta_2$ the \emph{loop parameters} of the planar algebra. In fact, any such planar algebra can be altered so $\delta_1=\delta_2=\delta$ as explained in \cite[Lemma 2.3.9]{palgnotes}. Then $Z_{\begin{tikzpicture}[scale=0.2]
		\draw[thick] (0,0) circle (1cm);
		\draw[shaded] (0,0) circle (0.5cm);
		\end{tikzpicture}}=\delta\cdot Z_{\begin{tikzpicture}[scale=0.1]
		\draw[thick] (0,0) circle (1cm);
		\end{tikzpicture}}$, which diagramatically can be expressed as:
	\[ \begin{tikzpicture}[scale=0.5,baseline=-0.4ex]
	\draw[thick] (0,0) circle (1cm);
	\draw[shaded] (0,0) circle (0.5cm);
	\end{tikzpicture} =\delta\cdot \begin{tikzpicture}[scale=0.5,baseline=-0.4ex]
	\draw[thick] (0,0) circle (1cm);
	\end{tikzpicture} \]
	We obtain a similar relation for the diagrams with the other shading.
	
	\begin{rem}
		The relation above allows us to eliminate loops from tangles. Using composition of tangles we can show that if $T$ is a tangle with a loop and $T'$ is the same tangle but without the loop, then $Z_T=\delta\cdot Z_{T'}$.
	\end{rem}
	
	For every $n\geq 0$ consider the following $n$-tangles:
	\[I^+_n=\begin{tikzpicture}[scale=0.5,baseline=0]
	\clip (0,0) circle (2.5cm);
	\draw[shaded] (-0.8,-3) rectangle (-0.4,3);
	\draw (0.8,-3) -- (0.8,3);
	\node at (0.25,1.5) {$\cdots$};
	\node at (0.25,-1.5) {$\cdots$};
	\node at (-180:2.2cm) {\scriptsize$\$$};
	\node at (-180:1.4cm) {\scriptsize$\$$};
	\node at (0,0) [Tbox, inner sep=4mm] {};
	\draw[very thick] (0,0) circle (2.48cm);
	\end{tikzpicture}\;,\quad I_n^-=\begin{tikzpicture}[scale=0.5,baseline=0,shaded]
	\clip (0,0) circle (2.5cm);
	\draw[shaded] (0,0) circle (2.48cm);
	\draw[fill=white] (-0.8,-3) rectangle (-0.4,3);
	\draw (0.8,-3) -- (0.8,3);
	\node at (0.25,1.5) {$\cdots$};
	\node at (0.25,-1.5) {$\cdots$};
	\node at (-180:2.2cm) {\scriptsize$\$$};
	\node at (-180:1.4cm) {\scriptsize$\$$};
	\node at (0,0) [Tbox, inner sep=4mm] {};
	\draw[very thick] (0,0) circle (2.48cm);
	\end{tikzpicture}
	 \]
	Clearly these tangles work as the identity under tangle composition. For example if $T$ is any positive $n$-tangle then $I_n^+\circ_1 T= T$ and consequently $Z_{I_n^+}\circ Z_T= Z_T$, where $Z_{I_n^+}:P_{n,+}\ra P_{n,+}$. This however does not imply that $Z_{I_n^+}$ has to be the identity map on $P_{n,+}$.
	
	\begin{defi}
		We say that a planar algebra is \emph{non-degenerate} if for all $n\geq 0$ we have that $Z_{I_n^\pm}=\ID_{P_{n,\pm}}$.
	\end{defi}
	
	It is also important to note that we have an algebra structure on every $P_{n,\pm}$, where the multiplication is determined by the multiplication tangles:
	
	\[m^+_k=\begin{tikzpicture}[scale=0.5,baseline=0]
	\clip (0,0) circle (2.5cm);
	\draw[shaded] (-0.7,-3) rectangle (-0.3,3);
	\draw (0.7,-3) -- (0.7,3);
	\node at (0.25,0) {$\cdots$};
	\node at (-180:2.2cm) {\scriptsize$\$$};
	\node at (-1.1,1.2) {\scriptsize$\$$};
	\node at (-1.1,-1.2) {\scriptsize$\$$};
	\node at (0,-1.2) [Tbox, inner sep=1mm] {$D_1$};
	\node at (0,1.2) [Tbox, inner sep=1mm] {$D_2$};
	\draw[very thick] (0,0) circle (2.48cm);
	\end{tikzpicture}\;,\quad m_k^-=\begin{tikzpicture}[scale=0.5,baseline=0,shaded]
	\clip (0,0) circle (2.5cm);
	\draw[shaded] (0,0) circle (2.48cm);
	\draw[fill=white] (-0.7,-3) rectangle (-0.3,3);
	\draw (0.7,-3) -- (0.7,3);
	\node at (0.25,0) {$\cdots$};
	\node at (-180:2.2cm) {\scriptsize$\$$};
	\node at (-1.1,1.2) {\scriptsize$\$$};
	\node at (-1.1,-1.2) {\scriptsize$\$$};
	\node at (0,-1.2) [Tbox, inner sep=3mm] {};
	\node at (0,1.2) [Tbox, inner sep=3mm] {};
	\draw[very thick] (0,0) circle (2.48cm);
	\end{tikzpicture}
	\]
	For $x_1,x_2\in P_{n,\pm}$ we set $x_1x_2:=Z_{m_n^\pm}(x_1\otimes x_2)$. The associativity of the product is obtained by composing the multiplication tangle with itself in the two possible ways and verifying that one obtains isotopic tangles. With this algebraic structure, we have that $P_{n,\pm}$ is a unital algebra for all $n\geq 0$ where the units are given by the distinguished elements defined by
	\[1^+_n=
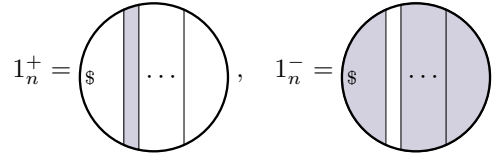
\begin{tikzpicture}[scale=0.5,baseline=0]
	\clip (0,0) circle (2cm);
	\draw[shaded] (-0.8,-3) rectangle (-0.4,3);
	\draw (0.8,-3) -- (0.8,3);
	\node at (0.25,0) {$\cdots$};
	\node at (-180:1.7cm) {\scriptsize$\$$};
	\draw[very thick] (0,0) circle (1.98cm);
	\end{tikzpicture}\;,\quad 1_n^-=\begin{tikzpicture}[scale=0.5,baseline=0]
	\clip (0,0) circle (2cm);
	\draw[shaded] (0,0) circle (1.98cm);
	\draw[fill=white] (-0.8,-3) rectangle (-0.4,3);
	\draw (0.8,-3) -- (0.8,3);
	\node at (0.25,0) {$\cdots$};
	\node at (-180:1.7cm) {\scriptsize$\$$};
	\draw[very thick] (0,0) circle (1.98cm);
	\end{tikzpicture}
	\]
	
	\begin{defi}[Planar $*$-algebra]\label{starpalg}
		We say that a planar algebra $P$ is a \emph{$*$-planar algebra} if each $P_{n,\pm}$ possesses a linear involution $*$ so that if $\theta$ is an orientation reversing diffeomorphism of $\R^2$ and $T$ a tangle with $k$ internal disks, then for any $x_1\otimes\cdots\otimes x_k\in \bigotimes_{i=1}^k P_{\partial(D_i)}$ we have
		\[ [Z_{\theta(T)}(x_1\otimes\cdots x_k)]^*=Z_T(x^*_1\otimes\cdots x^*_k). \]
	\end{defi}

        \begin{defi}[Jones projections]
		For $n\geq 1$, the \emph{Jones projections} in $P_{n+1,+}$ are $\delta^{-1} E_n$, where $E_n$ is the distinguished element associated to:
		\[ E_n=
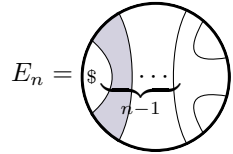
\begin{tikzpicture}[baseline=-0.4ex]
		\clip (0,0) circle (1cm); 
		\draw[shaded] (-150:1cm) .. controls ++(40:0.5cm) and ++(-40:0.5cm) .. (150:1cm) -- (150:1.5cm) -- (120:1cm) .. controls ++(-60:0.5cm) and ++(60:0.5cm) .. (-120:1cm) -- cycle;
		\node at (0,0) {$\cdots$};
		\node at (180:0.85cm) {\scriptsize$\$$};
		\draw (65:1cm) .. controls ++(-120:0.5cm) and ++(120:0.5cm) .. (-65:1cm);
		\draw (50:1cm) .. controls ++(-120:0.5cm) and ++(-170:0.5cm) .. (15:1cm);
		\draw (-50:1cm) .. controls ++(120:0.5cm) and ++(170:0.5cm) .. (-15:1cm);
		\node at (-0.2cm,-0.3cm) {$\underbrace{\hspace{3em}}_{n-1}$};
		\draw[very thick] (0,0) circle (0.98cm);
		\end{tikzpicture} \]
	\end{defi} 
 
	\begin{rem}
		In a $*$-planar algebra, $\delta^{-1}E_n$ is a projection in the classical sense.
	\end{rem}
	
	\begin{defi}[$C^*$-planar algebra]
		We say that a planar $*$-algebra is a \emph{$C^*$-planar algebra} if every $P_{n,\pm}$ has a norm which satisfies the $C^*$-identity with the algebra structure given by the multiplication tangle. That is, every $P_{n,\pm}$ is a $C^*$-algebra.  
	\end{defi}
	
	Since tangles only care about the isotopy class of the diagram, we can deform the disks to rectangles where half of the boundary points are at the top, the other half at the bottom and the marked interval is always on the left of the rectangle. Following this convention, composing two tangles with the multiplication tangle amounts to stacking them on top of each other. 
	
	\begin{rem}
		Using the convention from above, the involution from \ref{starpalg} can be diagramatically represented as
		
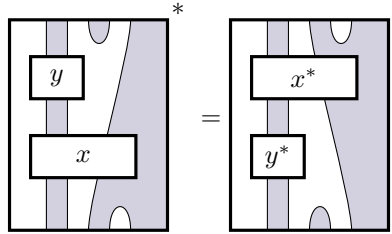
\[ \begin{tikzpicture}[scale=.7,baseline=1.4cm] 
		\def \wo {3}	
		\def \wi {1}	
		\def \ho {4}	
		\def \hi {1}	
		\def \wc {0.4}	
		\def \ws {0.4}	
		\def \wl {0.4+\ws}	
		\begin{scope}[shift={(0.7,0)}]							
		\draw[shaded] (0,0) -- ++(0,\ho) -- ++(\ws,0) -- ++(0,-\ho) -- 	cycle;
		\end{scope}	
		\begin{scope}[shift={(0.7+\wl,0)}]		
		\draw[shaded]  (0,0) .. controls ++(90:0.6) and ++(-90:0.6) .. (\wc+\wc,\ho) -- (\wo-1.5-\wl,\ho) -- (\wo-1.5-\wl,0) -- (\wc+\wc,0) .. controls ++(90:0.6) and ++(90:0.6) .. (\wc,0) -- cycle;
		\draw[shaded]  (\wc,\ho) .. controls ++(-90:0.6) and ++(-90:0.6) .. (0,\ho);
		\end{scope}	
		\node at (\wo+0.2,\ho+0.2) {$\ast$};
		\draw[very thick,fill=white] (0.4,2.5) rectangle ++(1,0.8);
		\draw[very thick,fill=white] (0.4,1) rectangle ++(2,0.8);
		\path (0.4,2.5) -- ++(1,0.8) node[midway] {$y$};
		\path (0.4,1) -- ++(2,0.8) node[midway] {$x$};
		\draw[very thick] (0,0) rectangle (\wo,\ho);
		\end{tikzpicture}= \begin{tikzpicture}[scale=.7,baseline=1.4cm] 
		\def \wo {3}	
		\def \wi {1}	
		\def \ho {4}	
		\def \hi {1}	
		\def \wc {0.4}	
		\def \ws {0.4}	
		\def \wl {0.4+\ws}	
		\begin{scope}[shift={(0.7,0)}]							
		\draw[shaded] (0,0) -- ++(0,\ho) -- ++(\ws,0) -- ++(0,-\ho) -- 	cycle;
		\end{scope}	
		\begin{scope}[shift={(0.7+\wl,0)}]		
		\draw[shaded]  (0,\ho) .. controls ++(-90:0.6) and ++(90:0.6) .. (\wc+\wc,0) -- (\wo-1.5-\wl,0) -- (\wo-1.5-\wl,\ho) -- (\wc+\wc,\ho) .. controls ++(-90:0.6) and ++(-90:0.6) .. (\wc,\ho) -- cycle;
		\draw[shaded]  (\wc,0) .. controls ++(90:0.6) and ++(90:0.6) .. (0,0);
		\end{scope}	
		\draw[very thick,fill=white] (0.4,1) rectangle ++(1,0.8);
		\draw[very thick,fill=white] (0.4,2.5) rectangle ++(2,0.8);
		\path (0.4,2.5) -- ++(2,0.8) node[midway] {$x^*$};
		\path (0.4,1) -- ++(1,0.8) node[midway] {$y^*$};
		\draw[very thick] (0,0) rectangle (\wo,\ho);
		\end{tikzpicture}\]
		which amounts to turning the tangle upside down and taking the involution of the elements on which we are evaluating.
	\end{rem}
	
	The planar algebra structure defined above implied that every $P_{n,\pm}$ is a unital associative algebra. There are however many tangles that induce maps between $P_{n,\pm}$ and $P_{n',\pm}$ where $n\neq n'$. One might wonder what extra structure these maps induce. There are two classes of tangles which begin to answer that question, the first being the \emph{inclusion tangles}:
	
	\[ \iota^+_{n}=\begin{tikzpicture}[scale=.7,baseline=4ex] 
	\def \wo {3.5}	
	\def \wi {2.4}	
	\def \ho {2}	
	\def \hi {1}	
	\def \wc {0.6}	
	\def \wl {\wi-0.6}	
	\def \ws {0.4}	
	\begin{scope}[shift={(0.2+\wi,0)}]		
	\draw  (\wc,0)--(\wc,\ho);
	\end{scope}	
	\begin{scope}[shift={(0.8,0)}]								
	\draw[shaded] (0,0) -- ++(0,\ho) -- ++(\ws,0) -- ++(0,-\ho) -- 	cycle;
	\draw (\wl,0) -- ++(0,\ho);
	\path (\ws+0.1,0.3) -- (\wl,0.3) node[midway] {$\cdots$};
	\path (\ws+0.1,\ho-0.25) -- (\wl,\ho-0.25) node[midway] {$\cdots$};
	\end{scope}	
	\draw[very thick] (0,0) rectangle (\wo,\ho);
	\draw[very thick,fill=white] (0.5,0.5) rectangle ++(\wi,\hi);
	\end{tikzpicture}\;, \quad  \iota^-_{n}=\begin{tikzpicture}[scale=.7,baseline=4ex] 
	\def \wo {3.5}	
	\def \wi {2.2}	
	\def \ho {2}	
	\def \hi {1}	
	\def \wc {0.6}	
	\def \wl {\wi}	
	\def \ws {0.6}	
	\begin{scope}[shift={(0.4,0)}]								
	\draw[shaded] (0,0) -- ++(0,\ho) -- ++(\ws,0) -- ++(0,-\ho) -- 	cycle;
	\draw (\wl,0) -- ++(0,\ho);
	\path (\ws+0.1,0.3) -- (\wl,0.3) node[midway] {$\cdots$};
	\path (\ws+0.1,\ho-0.25) -- (\wl,\ho-0.25) node[midway] {$\cdots$};
	\end{scope}	
	\draw[very thick] (0,0) rectangle (\wo,\ho);
	\draw[very thick,fill=white] (0.7,0.5) rectangle ++(\wi,\hi);
	\end{tikzpicture}\;,\quad \]
	where $\iota_{n}^{\pm}$ is a $n+1$-tangle. The other class of tangles are the \emph{conditional expectations}:
	\[ \alpha_n=\begin{tikzpicture}[scale=.7,baseline=4ex] 
	\def \wo {3.5}	
	\def \wi {2.4}	
	\def \ho {2}	
	\def \hi {1}	
	\def \wc {0.6}	
	\def \wl {\wi-0.9}	
	\def \ws {0.4}	
	\begin{scope}[shift={(0.2+\wi,0.5)}]		
	\draw (0,0) .. controls ++(-90:0.4) and ++(-90:0.4) .. (\wc,0)--(\wc,\hi) .. controls ++(90:0.4) and ++(90:0.4) .. (0,\hi) -- cycle;
	\end{scope}	
	\begin{scope}[shift={(0.8,0)}]								
	\draw[shaded] (0,0) -- ++(0,\ho) -- ++(\ws,0) -- ++(0,-\ho) -- 	cycle;
	\draw (\wl,0) -- ++(0,\ho);
	\path (\ws+0.1,0.3) -- (\wl,0.3) node[midway] {$\cdots$};
	\path (\ws+0.1,\ho-0.25) -- (\wl,\ho-0.25) node[midway] {$\cdots$};
	\end{scope}	
	\draw[very thick] (0,0) rectangle (\wo,\ho);
	\draw[very thick,fill=white] (0.5,0.5) rectangle ++(\wi,\hi);
	\end{tikzpicture}\;, \quad \beta_n=\begin{tikzpicture}[scale=.7,baseline=4ex] 
	\def \wo {3.5}	
	\def \wi {2.4}	
	\def \ho {2}	
	\def \hi {1}	
	\def \wc {0.6}	
	\def \wl {\wi-0.9}	
	\def \ws {1.2}	
	\begin{scope}[shift={(0,0)}]				
	\draw[shaded] (0,0) -- ++(0,\ho) -- ++(\ws,0)  -- ++(0,-\ho) -- 	cycle;
	\draw (\ws+\wl,0) -- ++(0,\ho);
	\path (\ws+0.1,0.3) -- ++(\wl,0) node[midway] {$\cdots$};
	\path (\ws+0.1,\ho-0.25) -- ++(\wl,0) node[midway] {$\cdots$};
	\end{scope}	
	\begin{scope}[shift={(0.3,0.5)}]		
	\draw[fill=white] (0,0) .. controls ++(-90:0.4) and ++(-90:0.4) .. (\wc,0)--(\wc,\hi) .. controls ++(90:0.4) and ++(90:0.4) .. (0,\hi) -- cycle;
	\end{scope}	
	\draw[very thick] (0,0) rectangle (\wo,\ho);
	\draw[very thick,fill=white] (0.6,0.5) rectangle ++(\wi,\hi);
	\end{tikzpicture}  \]
	where both are $n-1$-tangles. Note that these tangles only have one internal disk. Such tangles are called \emph{annular tangles}. When composing an annular tangle $T$ with any other tangle $S$ we will denote the composition simply by $T\circ S$ as there is no confusion as to where to place $S$. For such tangles $T$, $Z_T$ will simply be a linear map from $P_{\partial D}$ (where $D$ is the only internal disk of $T$) to $P_{\partial T}$. Therefore we can justify treating $T$ as a map $T: P_{\partial D}\ra P_{\partial T}$. With this convention, the tangles previously introduced define maps
	\begin{align*}
	\iota_{n}^+&:P_{n,+}\ra P_{n+1,+}\\
	\iota_{n}^-&:P_{n,-}\ra P_{n+1,+}\\
	\alpha_{n}&:P_{n,+}\ra P_{n-1,+}\\
	\beta_{n}&:P_{n,+}\ra P_{n-1,-}\\
	\end{align*}
	Observe that $\alpha_{n+1}\circ\iota_n^+=\delta I^+_n$ and $\beta_{n+1}\circ\iota_n^+=\delta I^-_n$. Therefore, if our planar algebra is non-degenerate, we have that $\iota_{n}^{\pm}$ is an injective map. In this case, we can treat $P_{n,+}$ as a subspace of $P_{n+1,+}$ and $P_{n,-}$ as a subspace of $P_{n+1,+}$, both invariant under the respective conditional expectation tangles.
	
	In principle, there are infinitely many tangles that be will inducing some structure on the subspaces $P_{n,\pm}$ and one might wonder if there is a nice set of tangles that generate all of them via composition. The next theorem, which can be found in \cite{kodiyalam2004jones}, answers that question.
	
	\begin{theo}\label{gentangles}
		Let $\T$ denote the set of all shaded tangles, and suppose $\T_1$ is a subclass of $\T$ which satisfies:
		\begin{enumerate}[(a)]
			\item $1^+_n,1^-_n\in \T_1$ for all $n\geq 0$,
			\item $E_n\in \T_1$ for all $n\geq 1$,
			\item $\alpha_n,\beta_n\in \T_1$ for all $n\geq 0$,
			\item $\iota^+_n,\iota^-_n\in \T_1$ for all $n\geq 0$,
			\item $m^+_n\in \T_1$ for all $n\geq 0$,
		\end{enumerate}
		and is closed under composition of tangles, then $\T=\T_1$.
	\end{theo}
	
	\begin{rem}
		This result will be extremely useful when verifying that something is a planar algebra. As long as we verify that the action of the generating tangles is compatible with composition, we can extend the action to all tangles and therefore obtain a planar algebra.
	\end{rem}

	Recall that $P_{0,\pm}$ can be identified with $\C$, hence $\alpha_1$ and $\beta_1$ define linear functionals. In fact, letting $\tr_n=\delta^{-n}\alpha_1\circ\alpha_2\circ \cdots \circ \alpha_n:P_{n,+}\ra P_{0,+}\simeq \C$ we obtain a linear functional on every space $P_{n,+}$. This linear functional is given by the following tangle:
	\[ \tr_n=\delta^{-n}\cdot\begin{tikzpicture}[baseline=1.2cm]
	\def \ch {0.8}
	\def \cw {0.4}
	\def \chh {0.8}
	\def \cww {1.4}
	\draw (1.4,0.7) .. controls ++(-90:0.4) and ++(-90:0.4) .. ++(\cw,0)-- ++(0,\ch) .. controls ++(90:0.4) and ++(90:0.4) .. ++(-\cw,0) -- cycle;
	\draw (0.8,0.7) .. controls ++(-90:0.8) and ++(-90:0.8) .. ++(\cww,0)-- ++(0,\chh) .. controls ++(90:0.8) and ++(90:0.8) .. ++(-\cww,0) -- cycle;
	\path (1,1.65) -- (1.3,1.65) node[midway] {$\cdots$};
	\path (1,0.55) -- (1.3,0.55) node[midway] {$\cdots$};
	\draw[very thick,fill=white] (0.6,0.7) rectangle ++(1,\ch);
	\draw[very thick] (0,0) rectangle (2.8,2.2);
	\end{tikzpicture} \]
	\begin{rem}
		It is clear that $\tr_n(xy)=\tr_n(yx)$ for all $x,y\in P_{n,+}$, as one can \emph{slide} $x$ across the strings from the bottom to the top. Moreover $\tr_n(1^+_n)=1$, hence $\tr_n$ is a normalized trace on $P_{n,+}$. However, if we draw a tangle as the one above but instead of capping the strings to the right we cap them to the left, we obtain another normalized trace which may not be equal to $\tr_n$.
	\end{rem}
	
	\begin{defi}[Spherical planar algebra]
		We say a planar algebra is \emph{spherical} if the two traces outlined in the previous remark coincide for all $n$.
	\end{defi}
	
	\begin{defi}[Positive planar algebra]
		Let $P$ be a C*-planar algebra, we say it is \emph{positive} if $\delta^n\tr_n:P_n\ra \C$ is a faithful positive map for all $n$.
	\end{defi}
	
	\begin{rem}
		The loop parameter of a positive planar algebra is positive.
	\end{rem}
	
	\begin{defi}
		Let $P_\bullet$ and $Q_\bullet$ be two planar algebras. A \emph{morphism} $\psi$ between $P_\bullet$ and $Q_\bullet$ is a family of linear maps $\psi_{n,\pm}:P_{n,\pm}\ra Q_{n,\pm}$ that is compatible with the action of tangles in the obvious way. We say $\psi:P_\bullet\ra Q_\bullet$ is an \emph{embedding} if every $\psi_{n,\pm}$ is injective.
	\end{defi}
	
	\subsection{Subfactor planar algebra}
	
	\begin{defi}
		A positive spherical planar algebra $P$ for which $\dim P_{n,\pm}<\infty$ and $\dim P_{0,\pm}=1$ is called a \emph{subfactor planar algebra}.
	\end{defi}

	Suppose $M_0\sst (M_1,\tr_1)$ is an extremal inclusion of \II\,factors, set $\delta=\sqrt{[M_1:M_0]}$ and let 
	\[ M_0\sst M_1\sst M_2 \sst \cdots \sst M_n\sst \cdots \]
	be the tower obtained from iterating the basic construction, this means $M_{n+1}=\langle M_{n},e_{n}\rangle$ where $e_n$ is the Jones projection with range in $L^2(M_{n-1},\tr_{n-1})$. Set $P_{n,+}=M_n\cap M_0'$ and $P_{n,-}=M_{n+1}\cap M'_1$ for $n\geq 0$. We want to endow $\{P_{n,\pm}\}$ with a planar algebra structure. As seen in \ref{gentangles}, to do this we just need to determine how the generating tangles act on these spaces. 
	
	\begin{enumerate}
		\item The distinguished element associated to $1_n^{\pm}$ is the unit $1$ in $M_n$.
		\item The distinguished element associated to $E_n$ is $\delta e_n\in P_{n+1,+}$.
		\item $m^+_n$ is given by the usual multiplication in $P_{n,+}$,
		\item $\alpha_n:P_{n,+}\ra P_{n-1,+}$ is given by $\delta E_{M_{n-1}}$ restricted to $P_{n,+}$, where $E_{M_{n-1}}:M_{n}\ra M_{n-1}$ is the $\tr_{M_n}$-preserving conditional expectation,
		\item $\beta_n:P_{n,+}\ra P_{n-1,-}$ is given by $\delta E_{M'_{1}}$ restricted to $P_{n,+}$, where $E_{M'_{1}}:M'_{0}\ra M'_{1}$ is the $\tr_{M'_0}$-preserving conditional expectation, where $\tr_{M'_0}|_{M_n}=\tr_{M_n}|_{M_0'}$ since $M_0\sst M_1$ is extremal,
		\item $\iota_n^{+}:P_{n,+}\ra P_{n+1,+}$ is given by the inclusion map $M_n \hookrightarrow M_{n+1}$ restricted to $P_{n,+}$,
		\item $\iota_n^{-}:P_{n,-}\ra P_{n+1,+}$ is given by the inclusion map $\iota:M'_1\hookrightarrow M'_0$ restricted to $P_{n,-}$.
	\end{enumerate}
	
	It is proved in \cite[Theorem 4.2.1]{jones1999planar} that the structure outlined above turns $\{P_{n,\pm}\}$ into a subfactor planar algebra. We denote this abstract planar algebra coming from the subfactor $M_0\sst M_1$ by $P^{M_0\sst M_1}$. A surprising fact is that all subfactor planar algebras will arise from an extremal subfactor in this way, as seen in \cite{popa1995axiomatization}, \cite{kodiyalam2009subfactor} or \cite{guionnet2010random}. 
	
	\subsection{Graph planar algebra}

	In this section we will focus on constructing a planar algebra from a bipartite graph. The construction for a general graph is similar but requires a more general definition of planar algebra which can be found in \cite{palgnotes}. It will not be needed here.
	
	We will be using the same notation used in Chapter 2. Let $\Gamma$ be a finite, connected, bipartite multigraph with adjacency matrix $\Lambda$.
	
	\begin{defi}
		If $T$ is a planar tangle, a \emph{state} $\sigma$ of $T$ is a function
		\[ \sigma: \{\text{regions of }T\}\cup \{\text{strings of T}\}\ra \V_+\cup\V_-\cup\E \]
		that sends unshaded regions to vertices in $\V_+$, shaded regions to vertices in $\V_-$ and strings to edges in $\E$ such that if $R_1$ and $R_2$ are two regions having a string $S$ as part of their boundary, then $\sigma(S)$ is an edge connecting $\sigma(R_1)$ and $\sigma(R_2)$.
	\end{defi}

	\begin{exa}
		Consider the following bipartite graph $\Gamma$,
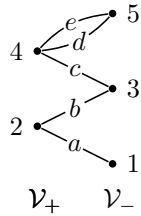
		\begin{center}
			\begin{tikzpicture}[semithick]
			\node[circle,fill=black,inner sep=0pt,minimum size=1mm,label=right:$3$] at (0,0) (A12) {};
			\node[circle,fill=black,inner sep=0pt,minimum size=1mm,label=left:$4$] at (-1,0.5) (A01)  {};
			\node[circle,fill=black,inner sep=0pt,minimum size=1mm,label=left:$2$] at (-1,-0.5) (A02) {};
			\node[circle,fill=black,inner sep=0pt,minimum size=1mm,label=right:$5$] at (0,1) (A11) {};
			\node[circle,fill=black,inner sep=0pt,minimum size=1mm,label=right:$1$] at (0,-1) (A13) {};
			\draw[postaction={decorate}] (A01) to [bend left =20] node[midway,circle,inner sep=0pt,fill=white] {$e$} (A11);
			\draw[postaction={decorate}] (A01) to [bend right =20] node[midway,circle,inner sep=0pt,fill=white] {$d$} (A11);
			\draw[postaction={decorate}] (A01) -- node[midway,circle,inner sep=0pt,fill=white] {$c$} (A12);
			\draw[postaction={decorate}] (A02) -- node[midway,circle,inner sep=0pt,fill=white] {$b$} (A12);
			\draw[postaction={decorate}] (A02) -- node[midway,circle,inner sep=0pt,fill=white] {$a$} (A13);
			\node at (-0.9,-1.5) {$\V_+$};
			\node at (-0.9,-1.5) {$\V_+$};
			\node at (0.1,-1.5) {$\V_-$};
			\end{tikzpicture}	
		\end{center}
		then we can represent a state $\sigma$ on a tangle $T$ in the following way:
		\[ \sigma=
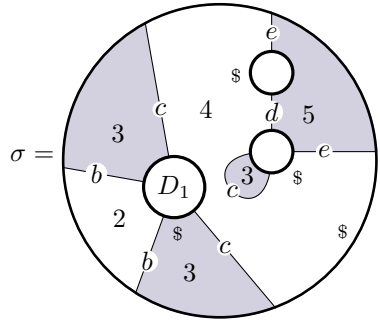
\begin{tikzpicture}[scale=.7,baseline=0]
		\clip (0,0) circle (3cm);
		
		\begin{scope}[shift=(10:1cm)] 
		\draw[shaded] (0,0)--node[midway,circle,inner sep=0pt,fill=white] {$e$} (0:2cm)--(45:3cm)--(90:3cm)--(0,0);  
		\draw[shaded] (0,0) .. controls ++(180:2cm) and ++(-90:2cm) .. node[midway,circle,inner sep=0pt,fill=white] {$c$} (0,0);
		\node at (45:1cm) {$5$};
		\node at (-135:0.65cm) {$3$};
		\node[circle,inner sep=0pt,fill=white] at (90:0.75cm) {$d$};
		\node[circle,inner sep=0pt,fill=white] at (90:2.25cm) {$e$};
		\end{scope}
		
		\begin{scope}[shift=(-150:1cm)]
		\draw[shaded] (0,0) -- (100:4cm) -- (170:4cm) -- (0,0);
		\draw[shaded] (0,0) -- (-110:4cm) -- (-50:4cm) -- (0,0);
		\node[circle,inner sep=0pt,fill=white] at (-110:1.5cm) {$b$};
		\node[circle,inner sep=0pt,fill=white] at (170:1.5cm) {$b$};
		\node[circle,inner sep=0pt,fill=white] at (100:1.5cm) {$c$};
		\node[circle,inner sep=0pt,fill=white] at (-50:1.5cm) {$c$};
		\end{scope}
		
		\node at (-105:2.2cm) {$3$};
		\node at (165:2cm) {$3$};
		\node at (210:2.2cm) {$2$};
		\node at (105:1cm) {$4$};
		
		\begin{scope}[shift=(10:1cm)] 
		\node at (0,0) [Tbox, inner sep=2mm] {};
		\node at (90:1.5cm) [Tbox, inner sep=2mm] {};
		\node at (-45:.7cm) {\scriptsize$\$$};
		\node at (115:1.6cm) {\scriptsize$\$$};
		\end{scope}
		\node at (-150:1cm) [Tbox, inner sep=1mm] {$D_1$};
		\node at (-120:1.6cm) {\scriptsize$\$$};
		\node at (-30:2.7cm) {\scriptsize$\$$};
		
		\draw[very thick] (0,0) circle (2.98cm);
		\end{tikzpicture} \]
	\end{exa}

	\begin{rem}
		For every disk in $T$, a state $\sigma$ on $T$ induces a pointed loop in $\Gamma$ by reading the values of the state clockwise around the disk starting from the marked interval. Reading around $D_1$ from the previous example we have:
		
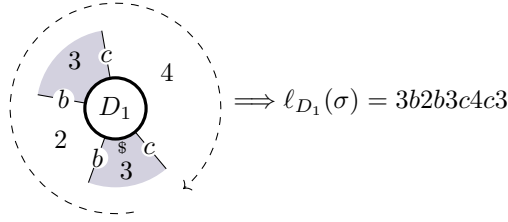
\[ \begin{tikzpicture}[scale=.7,baseline=0]
		\begin{scope}
			\clip (0,0) circle (1.5cm);
			\draw[shaded] (0,0) -- (100:4cm) -- (170:4cm) -- (0,0);
			\draw[shaded] (0,0) -- (-110:4cm) -- (-50:4cm) -- (0,0);
			\node[circle,inner sep=0pt,fill=white] at (-110:1cm) {$b$};
			\node[circle,inner sep=0pt,fill=white] at (170:1cm) {$b$};
			\node[circle,inner sep=0pt,fill=white] at (100:1cm) {$c$};
			\node[circle,inner sep=0pt,fill=white] at (-50:1cm) {$c$};
			\node at (0,0) [Tbox, inner sep=1mm] {$D_1$};
			\node at (-80:0.75cm) {\tiny$\$$};
			\node at (-80:1.2cm) {$3$};
			\node at (-150:1.2cm) {$2$};
			\node at (130:1.2cm) {$3$};
			\node at (35:1.2cm) {$4$};
		\end{scope}
		\draw[->,dashed] (-80:2cm) arc (280:-50:2cm);
		\end{tikzpicture} \Longrightarrow \ell_{D_1}(\sigma)=3b2b3c4c3 \]
		Given a disk $D$ in a tangle $T$, we denote by $\ell_D(\sigma)$ the loop associated to it, in particular we denote by $\ell_\sigma$ the loop associated to $D^T$.
	\end{rem}		
	
	Now let $\lambda:\V_+\cup \V_-\ra \R$ be defined as $\lambda|_{\V_+}=\lambda^0$ and $\lambda|_{\V_-}=d\lambda^1$ where $\lambda^0$ and $\lambda^1$ are the Markov trace vectors introduced in Definition \ref{markovTracev} and $d=\sqrt{\|\Lambda \Lambda^t\|}$.

	\begin{defi}\label{graphTaction}
		Given a state $\sigma$ on a tangle $T$ we define the \emph{rotation} $\Rot(\sigma)$ as follows:
		
		For every region $r$ in $T$ consider the positive orientation on it, then collapse every internal disk to a point so that each $r$ has as its boundary a union of oriented closed piecewise smooth curves. Let $\Rot(r)$ be the rotation number of the boundary of $r$. Then we set
			\[ \Rot(\sigma)=\prod_{\text{regions $r$ of }T}\lambda(\sigma(r))^{\Rot(r)} \]
	\end{defi}
	
	\begin{rem}
		The compute the rotation number of a region $r$, we orient $r$ positively and look at the orientation it induces on its boundary. Note that the boundary of $r$ is a finite collection of closed piecewise smooth curves. We define
		\[ \Rot(r)=\#\{\text{positively oriented curves}\} - \#\{\text{negatively oriented curves}\}.\]
		For example, if $r$ is
		
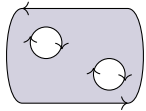
\[ \begin{tikzpicture}[scale=0.7]
			\draw[
			decoration={markings, mark=at position 0 with {\arrow{>}}, mark=at position 0.5 with {\arrow{>}}},
			postaction={decorate},shaded
			] (-45:1.2cm) .. controls ++(0:0.4) and ++(0:0.4) .. ++(0,1.8cm)-- ++(-2cm,0) .. controls ++(180:0.4) and ++(180:0.4) .. ++(0,-1.8cm) -- cycle;
			\begin{scope}[xshift=-0.7cm,yshift=0.3cm]
			\draw[
			decoration={markings, mark=at position 0 with {\arrow{<}}, mark=at position 0.5 with {\arrow{<}}},
			postaction={decorate},fill=white
			] (0,0) circle (0.3cm);
			\end{scope}
			\begin{scope}[xshift=0.5cm,yshift=-0.3cm]
			\draw[
			decoration={markings, mark=at position 0 with {\arrow{<}}, mark=at position 0.5 with {\arrow{<}}},
			postaction={decorate},fill=white
			] (0,0) circle (0.3cm);
			\end{scope}
		\end{tikzpicture} \]
		then $\Rot(r)=1-2=-1$.
	\end{rem}
	
	Recall that in Section 2.1, given a bipartite graph $\Gamma$, we defined the loop algebras $\{G_{n,\pm}\}$ whose elements are linear combinations of loops of length $2n$ starting in $\V_\pm$. Given $[\ell]\in G_{n,\pm}$ we define $\mu([\ell])=\lambda(v)\lambda(w)$ where $v$ and $w$ are the first (or last) and middle vertex of the loop $\ell$ and 
	\begin{align*}
	\lambda([\ell])&=\prod_{v\text{ vertex in }\ell}\lambda(v).
	\end{align*}
	To obtain the graph planar algebra we need to determine how planar tangles act on these spaces. Let $T$ be a tangle with $k$ internal disks $\{D_i\}_{i=1}^k$ and consider $[\ell_i]\in G_{\partial(D_i)}$ then we define
	\[ Z_T([\ell_1]\otimes \cdots \otimes [\ell_k])=\prod_i \sqrt{\frac{\mu([\ell_i])}{\lambda([\ell_i])}}\sum_{\substack{\sigma\text{ state on }T\\ \ell_i=\ell_{D_i(\sigma)}, \forall i}}\frac{\Rot(\sigma)}{\sqrt{\lambda([\ell_\sigma])\mu([\ell_{\sigma}])}} [\ell_\sigma].\]
	Therefore to evaluate $T$ we need consider all states $\sigma$ that induce precisely the loop $\ell_i$ around the disk $D_i$, then read the loop it induces around the output disk and multiply it by a correction factor. With this action one can show that $\{G_{n,\pm}\}$ has a planar algebra structure. The full proof can be found in \cite{palgnotes}.
	
	\begin{theo}
		Definition \ref{graphTaction} makes $\{G_{n,\pm}\}$ into a $C^*$-planar algebra.
	\end{theo}
	
	\begin{rem}
		Observe that for this planar algebra we do not necessarily have $\dim G_{0,\pm}\neq 1$, in fact $G_{0,\pm}\simeq \C^{|\V_\pm|}$. This means that if $|\V_\pm|>1$ then $\{G_{n,\pm}\}$ can not be a subfactor planar algebra as the notion of spherically does not make sense in this context.
	\end{rem}
	
	\begin{exa}
		Recall that tangles $T$ with no internal disks induce a distinguished element $Z_T(1)$. Hence, the tangle $Z_{\begin{tikzpicture}[scale=0.2]
			\draw[thick] (0,0) circle (1cm);
			\draw[shaded] (0,0) circle (0.5cm);
			\end{tikzpicture}}$ corresponds to an element in $G_{0,+}$. To compute this element we need to look at all possible states on this tangle. Let $\V_+=\{v_1^+,\dots,v_k^+\}$ and $\V_-=\{v_1^-,\dots,v_l^-\}$, every state in the tangle is given by a pair of vertices $(v_i^+,v_j^-)$ and an edge between them, hence:
		\[ Z_{\begin{tikzpicture}[scale=0.2]
			\draw[thick] (0,0) circle (1cm);
			\draw[shaded] (0,0) circle (0.5cm);
			\end{tikzpicture}}(1)=\sum_{i=1}^{k}\sum_{j=1}^{l}\sum_{\substack{\ve \text{ edge}\\ s(\ve)=v^+_i,\, t(\ve)=v_j^-}} \frac{d\lambda^1(v_j^-)}{\sqrt{\lambda^0(v_i^+)\lambda^0(v_i^+)}}[v^+_i].\]
		Here we are using the convention that $\mu([v_i^+])=\lambda(v_i^+)$. Recall that we have $\Lambda_{ij}$ edges from $v_i^+$ to $v_j^-$ and that $\Lambda \lambda^1=\lambda^0$. Therefore
		\begin{align*}
		Z_{\begin{tikzpicture}[scale=0.2]
			\draw[thick] (0,0) circle (1cm);
			\draw[shaded] (0,0) circle (0.5cm);
			\end{tikzpicture}}(1)&=\sum_{i=1}^{k}\sum_{j=1}^{l}d\Lambda_{ij}\frac{\lambda^1(v_j^-)}{\lambda^0(v_i^+)}[v^+_i]\\
		&=\sum_{i=1}^{k}d\frac{\lambda^0(v_i^+)}{\lambda^0(v_i^+)}[v^+_i]\\
		&=d\sum_{i=1}^{k}[v^+_i]=dZ_{\begin{tikzpicture}[scale=0.1]
			\draw[thick] (0,0) circle (1cm);
			\end{tikzpicture}}(1).
		\end{align*}
		Similarly we have $Z_{\begin{tikzpicture}[scale=0.2]
			\draw[thick,shaded] (0,0) circle (1cm);
			\draw[fill=white] (0,0) circle (0.5cm);
			\end{tikzpicture}}=d\cdot Z_{\begin{tikzpicture}[scale=0.1]
			\draw[thick,shaded] (0,0) circle (1cm);
			\end{tikzpicture}}$, so that our planar algebra has loop parameter $d$.
	\end{exa}
	
	\begin{exa}\label{gpaGenTang}
		We illustrate how some of the generating tangles act on this planar algebra. Consider the generating tangle
		\[ \alpha_{n}=
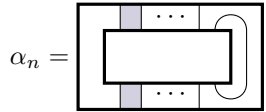
\begin{tikzpicture}[scale=.7,baseline=4ex] 
		\def \wo {3.5}	
		\def \wi {2.4}	
		\def \ho {2}	
		\def \hi {1}	
		\def \wc {0.6}	
		\def \wl {\wi-0.9}	
		\def \ws {0.4}	
		\begin{scope}[shift={(0.2+\wi,0.5)}]		
		\draw (0,0) .. controls ++(-90:0.4) and ++(-90:0.4) .. (\wc,0)--(\wc,\hi) .. controls ++(90:0.4) and ++(90:0.4) .. (0,\hi) -- cycle;
		\end{scope}	
		\begin{scope}[shift={(0.8,0)}]								
		\draw[shaded] (0,0) -- ++(0,\ho) -- ++(\ws,0) -- ++(0,-\ho) -- 	cycle;
		\draw (\wl,0) -- ++(0,\ho);
		\path (\ws+0.1,0.3) -- (\wl,0.3) node[midway] {$\cdots$};
		\path (\ws+0.1,\ho-0.25) -- (\wl,\ho-0.25) node[midway] {$\cdots$};
		\end{scope}	
		\draw[very thick] (0,0) rectangle (\wo,\ho);
		\draw[very thick,fill=white] (0.5,0.5) rectangle ++(\wi,\hi);
		\end{tikzpicture} \]
		acting on $[\ell]=[\ell_1 \ve_n^*\ve_{n+1}(\ell_2)^*]\in G_{n,+}$, where $n$ is even. Observe that if $\ve_{n}\neq \ve_{n+1}$ there exists no state on $\alpha_{n}$ which will induce the loop $\ell$ around the internal disk. On the other hand, if $\ve_{n}= \ve_{n+1}$ the only possible state $\sigma$ will induce the loop $[\ell_1(\ell_2)^{*}]$ on the external disk. Thus 
		\[ \alpha_{n}([\ell_1 \ve_n^*\ve_{n+1}(\ell_2)^*])=\delta_{\ve_n,\ve_{n+1}}\sqrt{\frac{\mu([\ell_1 \ve_n^*\ve_{n+1}(\ell_2)^*])}{\lambda([\ell_1 \ve_n^*\ve_{n+1}(\ell_2)^*])}}\frac{\Rot(\sigma)}{\sqrt{\lambda([\ell_1(\ell_2)^{*}])\mu([\ell_1(\ell_2)^{*}])}}[\ell_1(\ell_2)^{*}] \]
		Since every region in $\alpha_{n}$ has rotation number 1, we have that 
		\[ \Rot(\sigma)=\lambda([\ell_1(\ell_2)^*])\lambda(s(\ve_n)). \]
		Note that $\lambda([\ell_1 \ve_n^*\ve_{n+1}(\ell_2)^*])=\lambda([\ell_1(\ell_2)^{*}])\lambda(s(\ve_n))\lambda(t(\ve_{n}))$. Therefore 
		\[ \alpha_{n}([\ell_1 \ve_n^*\ve_{n+1}(\ell_2)^*])=\delta_{\ve_n,\ve_{n+1}}\sqrt{\frac{\mu([\ell_1 \ve_n^*\ve_{n+1}(\ell_2)^*])}{\mu([\ell_1(\ell_2)^{*}])}}\frac{\sqrt{\lambda(s(\ve_n))}}{\sqrt{\lambda(t(\ve_{n}))}}[\ell_1(\ell_2)^{*}]. \]
		Finally observe that the first vertex of $\ell_1 \ve_n^*\ve_{n+1}(\ell_2)^*$ and $\ell_1(\ell_2)^{*}$ are the same, whereas their middle vertices are $s(\ve_n)$ and $t(\ve_n)$, respectively. Hence we have 
		\[\alpha_{n}([\ell_1 \ve_n^*\ve_{n+1}(\ell_2)^*])=\delta_{\ve_n,\ve_{n+1}}\frac{\lambda(s(\ve_n))}{\lambda(t(\ve_{n}))}[\ell_1(\ell_2)^{*}]. \]
		Note that the formula above is very similar to the one we obtained in \ref{graphCExp}. In fact, we have 
		$$\alpha_{n}=\phi_{n-1,+}\circ(dE_{B_{n-1}})\circ (\phi_{n,+})^{-1}$$ 
		where $\phi_{n,+}$ is the isomorphism between $G_{n,+}$ and $B_0'\cap B_n$ described in Corollary \ref{loopIsom}.

        \end{exa}
        
		In a similar manner, we compute how the other generating tangles act on $G_{n,\pm}$:
		
		\begin{enumerate}
			\item $m^+_n$ will coincide with the multiplication in $G_{n,+}$,
			\item $\beta_n=(\phi_{n-1,-})\circ (dE_{B'_1})\circ (\phi_{n,+})^{-1}$, where $E_{B'_1}$ is the conditional expectation from $B'_0\cap B_n$ to $B'_1\cap B_n$ shown in \ref{loopcomCond},
			\item $\iota_n^+=(\phi_{n+1,+})\circ \iota \circ (\phi_{n,+})^{-1}$ where $\iota$ is the inclusion of $B_n$ into $B_{n+1}$,
			\item $\iota_n^-=(\phi_{n+1,+})\circ \iota \circ (\phi_{n,-})^{-1}$ where $\iota$ is the inclusion of $B'_1\cap B_{n+1}$ into $B'_0\cap B_n$.
			\item $E_n=(\phi_{n+1,+})(F_n)$.
		\end{enumerate}
	
	We deduce that the planar algebra structure on $\{G_{n,\pm}\}$ extends the structure we defined on Section 2.

\section{The embedding theorem}
Consider a nondegenerate commuting square of finite-dimensional 
C$^*$-algebras as in \ref{csfindim}: 
	\begin{equation}\label{csq}
	\begin{array}{ccc}
	A_{1,0} & \sst & A_{1,1} \\
	\cup &  & \cup \\
	A_{0,0} & \sst & A_{0,1}
	\end{array}
	\end{equation}
	
	Iterating the basic construction vertically and horizontally we obtain the following lattice of commuting squares:
	\[ \begin{array}{ccccccccc}
		A_{\infty,0} & \sst & A_{\infty,1} & \sst & A_{\infty,2} & \sst & \cdots & \sst & A_{\infty,\infty} \\
		\cup &  & \cup &  & \cup &  &  &  & \cup \\
		\vdots &  & \vdots &  & \vdots &  &  &  & \vdots \\
		\cup &  & \cup &  & \cup &  &  &  & \cup \\
		A_{2,0} & \sst & A_{2,1} & \sst & A_{2,2} & \sst & \cdots & \sst & A_{2,\infty}  \\
		\cup &  & \cup &  & \cup &  &  &  &  \\
		A_{1,0}& \sst & A_{1,1} & \sst & A_{1,2} & \sst & \cdots & \sst & A_{1,\infty} \\
		\cup &  & \cup &  & \cup &  &  &  & \cup \\
		A_{0,0} & \sst & A_{0,1} & \sst & A_{0,2} & \sst & \cdots & \sst & A_{0,\infty} 
	\end{array} \]
	where $A_{n,\infty}=\left(\bigcup_{k} A_{n,k}\right)''$ and $A_{\infty,k}=\left(\bigcup_{n} A_{n,k}\right)''$.
	
	In this section we will show that the subfactor planar algebra associated to $A_{0,\infty}\sst A_{1,\infty}$ embeds into the graph planar algebra associated to the Bratteli diagram $\Gamma$ of the first vertical inclusion 
$A_{0,0}\sst A_{1,0}$. A similar result holds of course for the vertical 
subfactor $A_{\infty, 0}\sst A_{\infty, 1}$. There are three important 
ideas used to prove this result. The first is that the relative commutants from the first vertical tower are isomorphic to the loop algebras associated to $\Gamma$ as shown in Section \ref{relcomloops}. The second idea is to use Ocneanu compactness (Theorem \ref{ocneanu}) 
	\[ A_{0,\infty}'\cap A_{n,\infty}=A_{0,1}'\cap A_{n,0}\sst A_{0,0}'\cap A_{n,0}. \]
	Consequently, the vector spaces of the subfactor planar algebra are included in the vector spaces of the graph planar algebra. The final idea is 
to use explicit computations for the actions of generating tangles for the graph planar algebra, like the ones in Example \ref{gpaGenTang}, to verify that these coincide with the action of the planar operad arising from the subfactor 
planar algebra.
	
	Let $P_\bullet=P^{A_{0,\infty}\sst A_{1,\infty}}$ be the subfactor planar algebra associated to $A_{0,\infty}\sst A_{1,\infty}$ and $P_{n,\pm}$, $n\geq 0$, be the corresponding vector spaces. Thus
	$P_{n,+}=A_{0,\infty}'\cap A_{n,\infty}$, and
	$P_{n,-}=A_{1,\infty}'\cap A_{n+1,\infty}$.

Moreover, the loop parameter for 
$P_\bullet$ is $\sqrt{[A_{1,\infty}:A_{0,\infty}]}$ .
	
	Now, consider $G_\bullet$ to be the graph planar algebra associated to $\Gamma$. Notice that the loop parameter is equal to $\sqrt{\|\Lambda\Lambda^t\|}$, where $\Lambda$ is the inclusion matrix for $A_{0,0}\sst A_{1,0}$. Since our commuting square is nondegenerate, we have
	\[ d^2=\|\Lambda\Lambda^t\|=\|\Lambda\|^2 =[A_{1,\infty}:A_{0,\infty}],\]  
	and therefore both planar algebras have the same loop parameter $d$.
	By Ocneanu compactness, we have 
	\[ A_{0,\infty}'\cap A_{n,\infty}\sst A_{0,0}'\cap A_{n,0}=Q_{n,+},\quad A_{1,\infty}'\cap  A_{n+1,\infty}\sst A_{1,0}'\cap A_{n+1,0}=Q_{n,-}.\]
	This means that the isomorphisms $\psi_{n,\pm}:Q_{n,\pm}\ra G_{n,\pm}$ from Remark \ref{towerisomor} define a map $\psi:P_\bullet\ra G_\bullet$. We will show that $\psi$ is a planar algebra embedding. Since $\psi_{n,\pm}$ is 
a $*$-algebra homomorphism for every $n$, we only need to verify that it 
preserves the actions of the following generating tangles:
	
	\begin{alignat*}{2}
	\alpha_n &= 
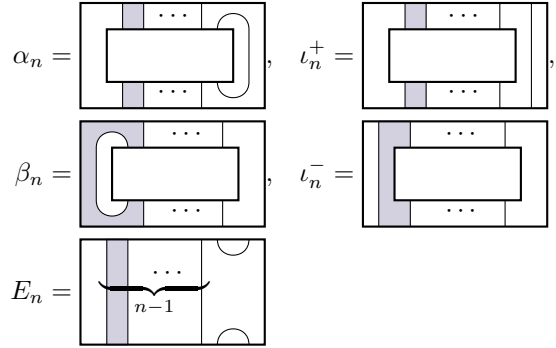
\begin{tikzpicture}[scale=.7,baseline=4ex] 
	\def \wo {3.5}	
	\def \wi {2.4}	
	\def \ho {2}	
	\def \hi {1}	
	\def \wc {0.6}	
	\def \wl {\wi-0.9}	
	\def \ws {0.4}	
	\begin{scope}[shift={(0.2+\wi,0.5)}]		
	\draw (0,0) .. controls ++(-90:0.4) and ++(-90:0.4) .. (\wc,0)--(\wc,\hi) .. controls ++(90:0.4) and ++(90:0.4) .. (0,\hi) -- cycle;
	\end{scope}	
	\begin{scope}[shift={(0.8,0)}]								
	\draw[shaded] (0,0) -- ++(0,\ho) -- ++(\ws,0) -- ++(0,-\ho) -- 	cycle;
	\draw (\wl,0) -- ++(0,\ho);
	\path (\ws+0.1,0.3) -- (\wl,0.3) node[midway] {$\cdots$};
	\path (\ws+0.1,\ho-0.25) -- (\wl,\ho-0.25) node[midway] {$\cdots$};
	\end{scope}	
	\draw[thick] (0,0) rectangle (\wo,\ho);
	\draw[thick,fill=white] (0.5,0.5) rectangle ++(\wi,\hi);
	\end{tikzpicture}, &\quad \iota_{n}^+&=\begin{tikzpicture}[scale=.7,baseline=4ex] 
	\def \wo {3.5}	
	\def \wi {2.4}	
	\def \ho {2}	
	\def \hi {1}	
	\def \wc {0.6}	
	\def \wl {\wi-0.6}	
	\def \ws {0.4}	
	\begin{scope}[shift={(0.2+\wi,0)}]		
	\draw  (\wc,0)--(\wc,\ho);
	\end{scope}	
	\begin{scope}[shift={(0.8,0)}]								
	\draw[shaded] (0,0) -- ++(0,\ho) -- ++(\ws,0) -- ++(0,-\ho) -- 	cycle;
	\draw (\wl,0) -- ++(0,\ho);
	\path (\ws+0.1,0.3) -- (\wl,0.3) node[midway] {$\cdots$};
	\path (\ws+0.1,\ho-0.25) -- (\wl,\ho-0.25) node[midway] {$\cdots$};
	\end{scope}	
	\draw[thick] (0,0) rectangle (\wo,\ho);
	\draw[thick,fill=white] (0.5,0.5) rectangle ++(\wi,\hi);
	\end{tikzpicture}, \\
	\beta_n&=\begin{tikzpicture}[scale=.7,baseline=4ex] 
	\def \wo {3.5}	
	\def \wi {2.4}	
	\def \ho {2}	
	\def \hi {1}	
	\def \wc {0.6}	
	\def \wl {\wi-0.9}	
	\def \ws {1.2}	
	\begin{scope}[shift={(0,0)}]				
	\draw[shaded] (0,0) -- ++(0,\ho) -- ++(\ws,0)  -- ++(0,-\ho) -- 	cycle;
	\draw (\ws+\wl,0) -- ++(0,\ho);
	\path (\ws+0.1,0.3) -- ++(\wl,0) node[midway] {$\cdots$};
	\path (\ws+0.1,\ho-0.25) -- ++(\wl,0) node[midway] {$\cdots$};
	\end{scope}	
	\begin{scope}[shift={(0.3,0.5)}]		
	\draw[fill=white] (0,0) .. controls ++(-90:0.4) and ++(-90:0.4) .. (\wc,0)--(\wc,\hi) .. controls ++(90:0.4) and ++(90:0.4) .. (0,\hi) -- cycle;
	\end{scope}	
	\draw[thick] (0,0) rectangle (\wo,\ho);
	\draw[thick,fill=white] (0.6,0.5) rectangle ++(\wi,\hi);
	\end{tikzpicture}, & \iota_n^{-}&=\begin{tikzpicture}[scale=.7,baseline=4ex] 
	\def \wo {3.5}	
	\def \wi {2.4}	
	\def \ho {2}	
	\def \hi {1}	
	\def \wc {0.6}	
	\def \wl {\wi-0.6}	
	\def \ws {0.6}	
	\begin{scope}[shift={(0.3,0)}]				
	\draw[shaded] (0,0) -- ++(0,\ho) -- ++(\ws,0)  -- ++(0,-\ho) -- 	cycle;
	\draw (\ws+\wl,0) -- ++(0,\ho);
	\path (\ws+0.1,0.3) -- ++(\wl,0) node[midway] {$\cdots$};
	\path (\ws+0.1,\ho-0.25) -- ++(\wl,0) node[midway] {$\cdots$};
	\end{scope}
	\draw[thick] (0,0) rectangle (\wo,\ho);
	\draw[thick,fill=white] (0.6,0.5) rectangle ++(\wi,\hi);
	\end{tikzpicture}\\
	E_n & =\begin{tikzpicture}[scale=.7,baseline=4ex] 
	\def \wo {3.5}	
	\def \wi {2.4}	
	\def \ho {2}	
	\def \hi {1}	
	\def \wc {0.6}	
	\def \wl {\wi-0.6}	
	\def \ws {0.4}	
	\begin{scope}[shift={(0.2+\wi,0)}]		
	\draw  (0,0) .. controls ++(90:0.4) and ++(90:0.4) .. (\wc,0);
	\draw  (\wc,\ho) .. controls ++(-90:0.4) and ++(-90:0.4) .. (0,\ho);
	\end{scope}	
	\begin{scope}[shift={(0.5,0)}]								
	\draw[shaded] (0,0) -- ++(0,\ho) -- ++(\ws,0) -- ++(0,-\ho) -- 	cycle;
	\draw (\wl,0) -- ++(0,\ho);
	\path (0,1) -- (\wl,1) node[midway] {$\underbrace{\hspace{2em}\cdots\hspace{1em}}_{n-1}$};
	\end{scope}	
	\draw[thick] (0,0) rectangle (\wo,\ho);
	\end{tikzpicture}
	\end{alignat*}

	Given a tangle $T$, we will use $Z^1_T$ to denote the action in $P_\bullet$ and $Z^2_T$ to denote the action in $G_\bullet$. In the case of the subfactor planar algebra $P_\bullet$, as stated in Section 3.1, we have:
	\begin{enumerate}[(a)]
		\item $Z^1_{\alpha_n}:P_{n,+}\ra P_{n-1,+}$ is given by $d E_{A_{n-1,\infty}}$ restricted to $P_{n,+}$, where  the $E_{A_{n-1,\infty}}$ is trace-preserving conditional expectation from $A_{n,\infty}$ to $A_{n-1,\infty}$.
		\item $Z^1_{\beta_n}:P_{n,+}\ra P_{n-1,-}$ is given by $d E_{A'_{1,\infty}}$, where $E_{A'_{1,\infty}}:A'_{0,\infty}\cap A_{n,\infty}\ra A'_{1,\infty}\cap A_{n,\infty}$ is the trace-preserving conditional expectation.
		\item $Z^1_{\iota_{n}^+}:P_{n,+}\ra P_{n+1,+}$ is given by the inclusion map from $A_{n,\infty}$ to $A_{n+1,\infty}$.
		\item $Z^1_{\iota_n^{-}}:P_{n,-}\ra P_{n+1,+}$ is given by the inclusion map from $A_{1,\infty}'\cap A_{n+1,\infty}$ to $A_{0,\infty}'\cap A_{n+1,\infty}$.
		\item $Z^1_{E_n}$ is given by the Jones projection in $A_{n+1,\infty}$ multiplied by the loop parameter $d$, that is $A_{n+1,\infty}=\langle A_{n,\infty}, E_n \rangle$
	\end{enumerate}
	From example \ref{gpaGenTang}, for the graph planar algebra $G_\bullet$ we have:
        \begin{enumerate}[(a)]
			\item $Z^2_{\alpha_{n}}=(\phi_{n-1,+})\circ(dE_{B_{n-1}})\circ (\phi_{n,+})^{-1}$,
			\item $Z^2_{\beta_n}=(\phi_{n-1,-})\circ (dE_{B'_1})\circ (\phi_{n,+})^{-1}$,
			\item $Z^2_{\iota_n^+}=(\phi_{n+1,+})\circ \iota \circ (\phi_{n,+})^{-1}$ where $\iota$ is the inclusion of $B_n$ into $B_{n+1}$,
			\item $Z^2_{\iota_n^-}=(\phi_{n+1,+})\circ \iota \circ (\phi_{n,-})^{-1}$ where $\iota$ is the inclusion of $B'_1\cap B_{n+1}$ into $B'_0\cap B_{n+1}$.
			\item $Z^2_{E_n}=(\phi_{n+1,+})(F_n)$ where $F_n$ is as in Definition \ref{JonesprojB}.
		\end{enumerate}
	Here $(B_n)_{n\geq 0}$ is the loop algebra associated to the first vertical inclusion as in Definition \ref{towerB}. We need to verify that all of these maps coincide when restricted to $P_\bullet$.
	
	It is clear that the inclusion map from $A_{n,\infty}$ to $A_{n+1,\infty}$ coincides with the inclusion map from $A_{n,0}$ to $A_{n+1,0}$ when restricted. From Proposition \ref{bcAandB}, we get that if $\iota$ is the inclusion from $B_n$ into $B_{n+1}$ then $(\varphi_{n+1})^{-1}\circ\iota\circ(\varphi_n)$  is the inclusion from $A_{n,0}$ into $A_{n+1,0}$ and therefore $(\varphi_{n+1})^{-1}\circ\iota\circ(\varphi_n)=Z^1_{\iota_{n}^+}:P_{n,+}\ra P_{n+1,+}$. From this we have
	\begin{align*}
	Z^2_{\iota_n^+}\circ \psi_{n,+} &= (\phi_{n+1,+})\circ \iota \circ (\phi_{n,+})^{-1}\circ \phi_{n,+}\circ \varphi_{n}\\
	&= (\phi_{n+1,+})\circ \iota\circ\varphi_{n}\\
	&= (\phi_{n+1,+})\circ\varphi_{n+1}\circ Z^1_{\iota_{n}^+}\\
	&=(\psi_{n+1,+})\circ Z^1_{\iota_{n}^+}
	\end{align*}
	which means $\iota_{n}^+$ acts the same with the subfactor planar algebra structure and the graph planar algebra structure. A similar computation shows that the actions by $i_n^{-}$ coincide.
	
	The commuting square relation implies $E_{A_{n-1,0}}(x)=E_{A_{n-1,\infty}}(x)$ for $x\in A_{n,0}$. Since the maps $\varphi_n$ are trace-preserving, we have
	$\varphi_{n-1}\circ E_{A_{n-1,0}}=E_{B_{n-1}}\circ \varphi_n$. This implies that for $x\in A_{0,\infty}'\cap A_{n,\infty}\sst A_{n,0}$, we have 
	\begin{align*}
	(Z^2_{\alpha_n}\circ \psi_{n,+})(x)&=(((\phi_{n-1,+})\circ(dE_{B_{n-1}})\circ (\phi_{n,+})^{-1})\circ (\phi_{n,+})\circ\varphi_{n})(x)\\
	&=((\phi_{n-1,+})\circ(dE_{B_{n-1}})\circ \varphi_{n})(x)\\
	&=((\phi_{n-1,+})\circ\varphi_{n-1}\circ dE_{A_{n-1,0}})(x)\\
	&=((\phi_{n-1,+})\circ\varphi_{n-1}\circ dE_{A_{n-1,\infty}})(x)\\
	&=(\psi_{n-1,+}\circ Z^1_{\alpha_n})(x)
	\end{align*}
	and therefore the actions by $\alpha_n$ coincide. Since the Jones projection that implements the conditional expectation from $A_{n,\infty}$ to $A_{n-1,\infty}$ is the same as the one that implements the conditional expectation from $A_{n,0}$ to $A_{n-1,0}$, a similar computation as the one above shows that the actions by $E_n$ also coincide. 
	
	Since our commuting square is nondegenerate, due to Corollary \ref{ppbasesCS}, we obtain a Pimsner-Popa basis $S$ that works for the inclusions $A_{0,0}\sst A_{1,0}$ and $A_{0,\infty}\sst A_{1,\infty}$. Using the formula from Proposition \ref{ppbasesrelcom}, for any $x\in A_{0,\infty}'\cap A_{n,\infty}\sst A'_{0,0}\cap A_{n,0}$ we have
	\[ E_{A_{1,\infty}'}(x)=\sum_{s\in S} sxs^*=E_{A_{1,0}'}(x). \]
	Observe that, since $\varphi_1(S)$ is also a Pimnser-Popa basis for $B_0\sst B_1$, for any $x\in A_{0,0}'\cap A_{n,0}$ we have
	\begin{align*}
	E_{B'_1}( \varphi_n(x))&= \sum_{s\in S} \varphi_1(s)\varphi_n(x)\varphi_1(s^*)\\
	&=\varphi_n\left(\sum_{s\in S} sxs^*\right)=\varphi_{n}\left(E_{A_{1,0}'}(x)\right).
	\end{align*}
	Therefore, if $x\in P_{n,+}=A_{0,\infty}'\cap A_{n,\infty}\sst A'_{0,0}\cap A_{n,0}$ we have 
	\begin{align*}
	(Z^2_{\beta_n}\circ \psi_{n,+})(x)&=\left((\phi_{n-1,-})\circ (dE_{B'_1})\circ (\phi_{n,+})^{-1}\circ (\phi_{n,+})\circ \varphi_{n}\right)(x)\\
	&=((\phi_{n-1,-})\circ(dE_{B'_1})\circ\varphi_{n})(x)\\
	&=((\phi_{n-1,-})\circ \varphi_{n}\circ dE_{A_{1,0}'})(x)\\
        &=((\psi_{n-1,-})\circ dE_{A_{1,0}'})(x)\\
	&=(\psi_{n-1,-})( dE_{A_{1,\infty}'}(x))= ((\psi_{n-1,-})\circ Z^{1}_{\beta_n})(x).
	\end{align*}
	From this, we conclude that the action by $\beta_n$ on $P_{n,+}$ is 
the same for the subfactor planar algebra structure or the graph planar 
algebra structure. Thus we have proved:
	\begin{theo}\label{embeddingT}
	Consider commuting square as in (\ref{csq}). Let $P_\bullet$ be 
the subfactor planar algebra associated to $A_{0,\infty}\sst A_{1,\infty}$ 
and $G_\bullet$ the graph planar algebra associated to the Bratteli 
diagram of $A_{0,0}\sst A_{1,0}$. The map $\psi:P_\bullet\ra G_\bullet$ 
defined above is a planar algebra embedding.
	\end{theo}
A similar result holds of course for the vertical subfactor
$A_{\infty,0}\sst A_{\infty,1}$.

\section{A hyperfinite $A_\infty$ subfactor with index $\frac{5+\sqrt{13}}{2}$}

In this section, we will apply theorem \ref{embeddingT} to 
show the existence of an irreducible, hyperfinite subfactor with 
Temperley-Lieb-Jones (also called ``trivial'' or ``minimal''))
standard invariant and Jones index equal to that of the 
Haagerup subfactor \cite{asaedahaagerup1999}. Recall that the
Haagerup subfactor is special in that it is the finite depth subfactor 
with smallest Jones index $>4$. It has index $\frac{5+\sqrt{13}}{2}$,
and all irreducible
subfactors with index between $4$ and $\frac{5+\sqrt{13}}{2}$ have
trivial standard invariant, and hence have $A_\infty$ principal graphs. 
Note that Popa showed in 
\cite{popa1993markov} that for any $\lambda>4$, one can construct 
irreducible \emph{non-hyperfinite} subfactors with index $\lambda$ and 
trivial standard invariant. Whether this can be done in the hyperfinite
world remains a major open problem. 

We establish here that this is possible
for $\lambda = \frac{5+\sqrt{13}}{2}$ and show in future work that
the result also holds for all finite depth Jones indices between
$4$ and $5$, and for $3 + \sqrt{5}$ \cite{bischcaceres2}. In fact, 
we conjecture that every index of an irreducible, hyperfinite subfactor 
can be realized by one with trivial standard invariant. Our evidence 
for this conjecture
is a little thin, as we are only able to establish it for finite depth,
small index subfactors. We strongly believe it should hold at least for 
Jones indices of finite depth subfactors.


\begin{theo} There is a subfactor of the hyperfinite II$_1$ factor with
index $\frac{5+\sqrt{13}}{2}$ and
Temperley-Lieb-Jones standard invariant, and hence A$_\infty$ principal
graphs.
\end{theo}
\begin{proof}
Haagerup and Schou constructed a hyperfinite subfactor 
$N\sst M$ with index $\frac{5+\sqrt{13}}{2}$ from a nondegenerate 
commuting square of
multi-matrix algebras \cite[Chapter 7]{schou2013commuting} whose 
first vertical inclusion graph is given by
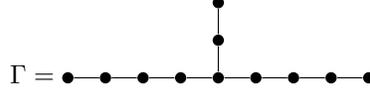
\[ \Gamma=\begin{tikzpicture}[scale=0.5]
	\node[circle,fill=black,inner sep=0pt,minimum size=1.5mm] at (1,0) (A2) {};
	\node[circle,fill=black,inner sep=0pt,minimum size=1.5mm] at (2,0) (A3) {};
	\node[circle,fill=black,inner sep=0pt,minimum size=1.5mm] at (3,0) (A4) {};
	\node[circle,fill=black,inner sep=0pt,minimum size=1.5mm] at (4,0) (A5) {};
	\node[circle,fill=black,inner sep=0pt,minimum size=1.5mm] at (5,1) (B1) {};
	\node[circle,fill=black,inner sep=0pt,minimum size=1.5mm] at (5,2) (B2) {};
	\node[circle,fill=black,inner sep=0pt,minimum size=1.5mm] at (5,0) (P) {};
	\node[circle,fill=black,inner sep=0pt,minimum size=1.5mm] at (6,0) (C1) {};
	\node[circle,fill=black,inner sep=0pt,minimum size=1.5mm] at (7,0) (C2) {};
	\node[circle,fill=black,inner sep=0pt,minimum size=1.5mm] at (8,0) (C3) {};
	\node[circle,fill=black,inner sep=0pt,minimum size=1.5mm] at (9,0) (C4) {};
	\draw (A2) -- (A3) -- (A4) -- (A5) -- (P) -- (C1) -- (C2) -- (C3) -- (C4);
	\draw (P) -- (B1) -- (B2);
\end{tikzpicture} \] 

This commuting square subfactor is irreducible by Wenzl's criterion
\cite{wenzl1988hecke}.

By the embedding theorem \ref{embeddingT}, the subfactor planar algebra 
of $N\sst M$ must embed
into the graph planar algebra of $\Gamma$. Since Peters showed 
that the planar algebra of the Haagerup subfactor does not embed 
in the graph planar algebra of $\Gamma$ (Theorem 6.8 in 
\cite{peters2010planar}), we can conclude that the
commuting square subfactor $N\sst M$ cannot be the Haagerup subfactor 
and therefore, by the classification of small index subfactor planar 
algebras, must be an $A_\infty$ subfactor.

Another way to see that the Haagerup subfactor planar algebra does not
embed in the graph planar algebra of $\Gamma$ uses corollary 1.4 
in \cite{grossman2018extended}, which proves that the Haagerup 
subfactor planar algebra only embeds in the graph planar algebra associated 
to one of the following three graphs:

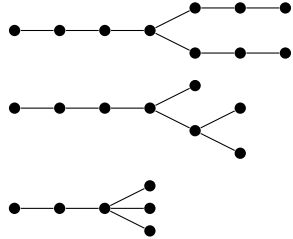
\begin{align*}
& \begin{tikzpicture}[scale=0.6]
\node[circle,fill=black,inner sep=0pt,minimum size=1.5mm] at (1,0) (A1) {};
\node[circle,fill=black,inner sep=0pt,minimum size=1.5mm] at (2,0) (A2) {};
\node[circle,fill=black,inner sep=0pt,minimum size=1.5mm] at (3,0) (A3) {};
\node[circle,fill=black,inner sep=0pt,minimum size=1.5mm] at (5,0.5) (B1) {};
\node[circle,fill=black,inner sep=0pt,minimum size=1.5mm] at (6,0.5) (B2) {};
\node[circle,fill=black,inner sep=0pt,minimum size=1.5mm] at (7,0.5) (B3) {};
\node[circle,fill=black,inner sep=0pt,minimum size=1.5mm] at (4,0) (P) {};
\node[circle,fill=black,inner sep=0pt,minimum size=1.5mm] at (5,-0.5) (C1) {};
\node[circle,fill=black,inner sep=0pt,minimum size=1.5mm] at (6,-0.5) (C2) {};
\node[circle,fill=black,inner sep=0pt,minimum size=1.5mm] at (7,-0.5) (C3) {};
\draw (A1) -- (A2) -- (A3) -- (P) -- (B1) -- (B2) -- (B3);
\draw (P) -- (C1) -- (C2) -- (C3);
\end{tikzpicture}\\
& \begin{tikzpicture}[scale=0.6]
\node[circle,fill=black,inner sep=0pt,minimum size=1.5mm] at (1,0) (A1) {};
\node[circle,fill=black,inner sep=0pt,minimum size=1.5mm] at (2,0) (A2) {};
\node[circle,fill=black,inner sep=0pt,minimum size=1.5mm] at (3,0) (A3) {};
\node[circle,fill=black,inner sep=0pt,minimum size=1.5mm] at (5,0.5) (B1) {};
\node[circle,fill=black,inner sep=0pt,minimum size=1.5mm] at (4,0) (P) {};
\node[circle,fill=black,inner sep=0pt,minimum size=1.5mm] at (5,-0.5) (C1) {};
\node[circle,fill=black,inner sep=0pt,minimum size=1.5mm] at (6,0) (C2) {};
\node[circle,fill=black,inner sep=0pt,minimum size=1.5mm] at (6,-1) (C3) {};
\draw (A1) -- (A2) -- (A3) -- (P) -- (B1);
\draw (P) -- (C1) -- (C2);
\draw (C1) -- (C3);
\end{tikzpicture}\\
& \begin{tikzpicture}[scale=0.6]
\node[circle,fill=black,inner sep=0pt,minimum size=1.5mm] at (2,0) (A2) {};
\node[circle,fill=black,inner sep=0pt,minimum size=1.5mm] at (3,0) (A3) {};
\node[circle,fill=black,inner sep=0pt,minimum size=1.5mm] at (5,0.5) (B1) {};
\node[circle,fill=black,inner sep=0pt,minimum size=1.5mm] at (4,0) (P) {};
\node[circle,fill=black,inner sep=0pt,minimum size=1.5mm] at (5,-0.5) (C1) {};
\node[circle,fill=black,inner sep=0pt,minimum size=1.5mm] at (5,0) (C2) {};
\draw (A2) -- (A3) -- (P) -- (B1);
\draw (P) -- (C1);
\draw (P) -- (C2);
\end{tikzpicture}
\end{align*}

Note that the first two graphs are the two principal graphs of the Haagerup 
subfactor. Since $\Gamma$ is not on this list, we are done.
\end{proof}

This is not the only index where this technique can be applied. From 
theorem 1.3 in \cite{grossman2018extended} we know that the Extended 
Haagerup subfactor planar algebra can only embed in the graph planar 
algebras of the following graphs:
\begin{table}[h!]
\centering
\begin{tblr}[caption={Extended Haagerup module graphs}]{
colspec={lrc},
vline{3},
hline{2},
row{1,2}={abovesep=5pt,belowsep=5pt}
}
Principal graph&\begin{tikzpicture}[scale=0.5,baseline=0]
\node[circle,fill=black,inner sep=0pt,minimum size=1.5mm] at (-1,0) (D1) {};
\node[circle,fill=black,inner sep=0pt,minimum size=1.5mm] at (-2,0) (D2) {};
\node[circle,fill=black,inner sep=0pt,minimum size=1.5mm] at (-3,0) (D3) {};
\node[circle,fill=black,inner sep=0pt,minimum size=1.5mm] at (0,0) (A0) {};
\node[circle,fill=black,inner sep=0pt,minimum size=1.5mm] at (1,0) (A1) {};
\node[circle,fill=black,inner sep=0pt,minimum size=1.5mm] at (2,0) (A2) {};
\node[circle,fill=black,inner sep=0pt,minimum size=1.5mm] at (3,0) (A3) {};
\node[circle,fill=black,inner sep=0pt,minimum size=1.5mm] at (5,0.5) (B1) {};
\node[circle,fill=black,inner sep=0pt,minimum size=1.5mm] at (6,0.5) (B2) {};
\node[circle,fill=black,inner sep=0pt,minimum size=1.5mm] at (7,0.5) (B3) {};
\node[circle,fill=black,inner sep=0pt,minimum size=1.5mm] at (4,0) (P) {};
\node[circle,fill=black,inner sep=0pt,minimum size=1.5mm] at (5,-0.5) (C1) {};
\node[circle,fill=black,inner sep=0pt,minimum size=1.5mm] at (6,-0.5) (C2) {};
\node[circle,fill=black,inner sep=0pt,minimum size=1.5mm] at (7,-0.5) (C3) {};
\draw (D3) -- (D2) -- (D1) -- (A0) -- (A1) -- (A2) -- (A3) -- (P) -- (B1) -- (B2) -- (B3);
\draw (P) -- (C1) -- (C2) -- (C3);
\end{tikzpicture} &\begin{tikzpicture}[scale=0.5,baseline=0]
\node[circle,fill=black,inner sep=0pt,minimum size=1.5mm] at (1,1) (A1) {};
\node[circle,fill=black,inner sep=0pt,minimum size=1.5mm] at (2,0.5) (A2) {};
\node[circle,fill=black,inner sep=0pt,minimum size=1.5mm] at (1,-1) (B1) {};
\node[circle,fill=black,inner sep=0pt,minimum size=1.5mm] at (2,-0.5) (B2) {};
\node[circle,fill=black,inner sep=0pt,minimum size=1.5mm] at (3,0) (C1) {};
\node[circle,fill=black,inner sep=0pt,minimum size=1.5mm] at (4,0) (C2) {};
\node[circle,fill=black,inner sep=0pt,minimum size=1.5mm] at (5,0) (C3) {};
\node[circle,fill=black,inner sep=0pt,minimum size=1.5mm] at (6,0) (C4) {};
\node[circle,fill=black,inner sep=0pt,minimum size=1.5mm] at (7,0.5) (D1) {};
\node[circle,fill=black,inner sep=0pt,minimum size=1.5mm] at (7,-0.5) (D2) {};
\node[circle,fill=black,inner sep=0pt,minimum size=1.5mm] at (8,-1) (D3) {};
\draw (A1) -- (A2) -- (C1) -- (C2) -- (C3) -- (C4) -- (D1);
\draw (B1) -- (B2) -- (C1);
\draw (C4) -- (D2) -- (D3);
\end{tikzpicture}\\
Dual principal graph&\begin{tikzpicture}[scale=0.5,baseline=0]
\node[circle,fill=black,inner sep=0pt,minimum size=1.5mm] at (-1,0) (D1) {};
\node[circle,fill=black,inner sep=0pt,minimum size=1.5mm] at (-2,0) (D2) {};
\node[circle,fill=black,inner sep=0pt,minimum size=1.5mm] at (-3,0) (D3) {};
\node[circle,fill=black,inner sep=0pt,minimum size=1.5mm] at (0,0) (A0) {};
\node[circle,fill=black,inner sep=0pt,minimum size=1.5mm] at (1,0) (A1) {};
\node[circle,fill=black,inner sep=0pt,minimum size=1.5mm] at (2,0) (A2) {};
\node[circle,fill=black,inner sep=0pt,minimum size=1.5mm] at (3,0) (A3) {};
\node[circle,fill=black,inner sep=0pt,minimum size=1.5mm] at (5,0.5) (B1) {};
\node[circle,fill=black,inner sep=0pt,minimum size=1.5mm] at (4,0) (P) {};
\node[circle,fill=black,inner sep=0pt,minimum size=1.5mm] at (5,-0.5) (C1) {};
\node[circle,fill=black,inner sep=0pt,minimum size=1.5mm] at (6,0) (C2) {};
\node[circle,fill=black,inner sep=0pt,minimum size=1.5mm] at (6,-1) (C3) {};
\draw (D3) -- (D2) -- (D1) -- (A0) -- (A1) -- (A2) -- (A3) -- (P) -- (B1);
\draw (P) -- (C1) -- (C2);
\draw (C1) -- (C3);
\end{tikzpicture}&\begin{tikzpicture}[scale=0.5,baseline=0]
\node[circle,fill=black,inner sep=0pt,minimum size=1.5mm] at (2,0.5) (A2) {};
\node[circle,fill=black,inner sep=0pt,minimum size=1.5mm] at (2,-0.5) (B2) {};
\node[circle,fill=black,inner sep=0pt,minimum size=1.5mm] at (3,0) (C1) {};
\node[circle,fill=black,inner sep=0pt,minimum size=1.5mm] at (4,0) (C2) {};
\node[circle,fill=black,inner sep=0pt,minimum size=1.5mm] at (5,0) (C3) {};
\node[circle,fill=black,inner sep=0pt,minimum size=1.5mm] at (6,0) (C4) {};
\node[circle,fill=black,inner sep=0pt,minimum size=1.5mm] at (7,0) (C5) {};
\node[circle,fill=black,inner sep=0pt,minimum size=1.5mm] at (8,0) (C6) {};
\node[circle,fill=black,inner sep=0pt,minimum size=1.5mm] at (9,0.5) (D1) {};
\node[circle,fill=black,inner sep=0pt,minimum size=1.5mm] at (9,-0.5) (D2) {};
\node[circle,fill=black,inner sep=0pt,minimum size=1.5mm] at (6,-1) (E1) {};
\draw (A2) -- (C1) -- (C2) -- (C3) -- (C4) -- (C5) -- (C6) -- (D1);
\draw (B2) -- (C1);
\draw (C6) -- (D2);
\draw (C4) -- (E1);
\end{tikzpicture}
\end{tblr}
\caption{Extended Haagerup module graphs}
\label{exthaagMGs}
\end{table}

More generally, Theorem 1.2 in \cite{grossman2018extended} states:
\begin{theo}\label{fusionembed}
        Suppose $P_\bullet$ is a finite depth subfactor planar algebra. Let $\mathcal{C}$ denote the unitary multifusion category of projections in $P_\bullet$, with distinguished object $X=\ID_{1,+}\in P_{1,+}$, and the standard unitary pivotal structure with respect to $X$. There is an equivalence between:
        \begin{enumerate}
            \item Planar algebra embeddings $P_\bullet \ra G_\bullet$, where $G_\bullet$ is the graph planar algebra associated to a finite connected bipartite graph $\Gamma$, and
            \item indecomposable finitely semisimple pivotal left $\mathcal{C}$-module $C^*$ categories $\mathcal{M}$ whose fusion graph with respect to $X$ is $\Gamma$.
        \end{enumerate}
    \end{theo}

In a forthcoming article, we have constructed commuting squares whose  
Bratteli diagram for the first vertical inclusion has norm squared
$\sim 4.37720$ ($=$ the index of the Extended Haagerup subfactor),
resp. $\frac{5+\sqrt{17}}{2}$ ($=$ the index of the Asaeda-Haagerup
subfactor) and is not one of 
the module graphs as in (2) of the previous theorem for these two
subfactors. This establishes then the existence of new irreducible, 
hyperfinite subfactors at these indices, with trivial standard invariant.


\nocite{*}
\bibliographystyle{alpha}
\bibliography{references}{}	

\newcommand{\etalchar}[1]{$^{#1}$}
\begin{thebibliography}{GMP{\etalchar{+}}23}

\bibitem[AH99]{asaedahaagerup1999}
M.~Asaeda and U.~Haagerup.
\newblock Exotic subfactors of finite depth with {J}ones indices
  {$(5+\sqrt{13})/2$} and {$(5+\sqrt{17})/2$}.
\newblock {\em Comm. Math. Phys.}, 202(1):1--63, 1999.

\bibitem[AMP23]{afzaly2015classification}
Narjess Afzaly, Scott Morrison, and David Penneys.
\newblock The classification of subfactors with index at most {$5\frac14$}.
\newblock {\em Mem. Amer. Math. Soc.}, 284(1405):v+81, 2023.

\bibitem[BC]{bischcaceres2}
Dietmar Bisch and Julio C{\'a}ceres.
\newblock New hyperfinite subfactors with infinite depth.
\newblock in preparation.

\bibitem[Bis94]{bisch1994}
Dietmar Bisch.
\newblock An example of an irreducible subfactor of the hyperfinite {${\rm
  II}_1$} factor with rational, noninteger index.
\newblock {\em J. Reine Angew. Math.}, 455:21--34, 1994.

\bibitem[Bis97]{bisch1997bimodules}
Dietmar Bisch.
\newblock Bimodules, higher relative commutants and the fusion algebra
  associated to a subfactor.
\newblock In {\em Operator algebras and their applications ({W}aterloo, {ON},
  1994/1995)}, volume~13 of {\em Fields Inst. Commun.}, pages 13--63. Amer.
  Math. Soc., Providence, RI, 1997.

\bibitem[BJ97]{bisch1997intermediate}
Dietmar Bisch and Vaughan Jones.
\newblock Algebras associated to intermediate subfactors.
\newblock {\em Invent. Math.}, 128(1):89--157, 1997.

\bibitem[BPMS12]{bigelow2012extendedHaagerup}
Stephen Bigelow, Emily Peters, Scott Morrison, and Noah Snyder.
\newblock {Constructing the extended Haagerup planar algebra}.
\newblock {\em Acta Mathematica}, 209(1):29 -- 82, 2012.

\bibitem[Bur15]{burstein2015}
Richard~D. Burstein.
\newblock Group-type subfactors and {H}adamard matrices.
\newblock {\em Trans. Amer. Math. Soc.}, 367(10):6783--6807, 2015.

\bibitem[Che24]{chen2024stardardlattices}
Quan Chen.
\newblock Standard-lattices, rigid {$C^*$}-tensor categories, and (bi)modules.
\newblock {\em Documenta Matematica}, 29(2):247--341, 2024.

\bibitem[CHPS21]{coleshuston2021}
Desmond Coles, Peter Huston, David Penneys, and Srivatsa Srinivas.
\newblock The module embedding theorem via towers of algebras.
\newblock {\em J. Funct. Anal.}, 280(11):Paper No. 108965, 52, 2021.

\bibitem[EGNO15]{etingof2016tensor}
Pavel Etingof, Shlomo Gelaki, Dmitri Nikshych, and Victor Ostrik.
\newblock {\em Tensor categories}, volume 205 of {\em Mathematical Surveys and
  Monographs}.
\newblock American Mathematical Society, Providence, RI, 2015.

\bibitem[EK98]{evanskawahigashi1998}
David~E. Evans and Yasuyuki Kawahigashi.
\newblock {\em Quantum symmetries on operator algebras}.
\newblock Oxford Mathematical Monographs. The Clarendon Press, Oxford
  University Press, New York, 1998.
\newblock Oxford Science Publications.

\bibitem[ENO05]{etingof2005fusion}
Pavel Etingof, Dmitri Nikshych, and Viktor Ostrik.
\newblock On fusion categories.
\newblock {\em Ann. of Math. (2)}, 162(2):581--642, 2005.

\bibitem[ENO10]{etingof2010fusion}
Pavel Etingof, Dmitri Nikshych, and Victor Ostrik.
\newblock Fusion categories and homotopy theory.
\newblock {\em Quantum Topol.}, 1(3):209--273, 2010.
\newblock With an appendix by Ehud Meir.

\bibitem[GdlHJ89]{goodman2012coxeter}
Frederick~M. Goodman, Pierre de~la Harpe, and Vaughan Jones.
\newblock {\em Coxeter graphs and towers of algebras}, volume~14 of {\em
  Mathematical Sciences Research Institute Publications}.
\newblock Springer-Verlag, New York, 1989.

\bibitem[GIS18]{grossman2018asaeda}
Pinhas Grossman, Masaki Izumi, and Noah Snyder.
\newblock The {A}saeda-{H}aagerup fusion categories.
\newblock {\em J. Reine Angew. Math.}, 743:261--305, 2018.

\bibitem[GJS10]{guionnet2010random}
Alice Guionnet, Vaughan. Jones, and Dimitri Shlyakhtenko.
\newblock Random matrices, free probability, planar algebras and subfactors.
\newblock In {\em Quanta of maths}, volume~11 of {\em Clay Math. Proc.}, pages
  201--239. Amer. Math. Soc., Providence, RI, 2010.

\bibitem[GMP{\etalchar{+}}23]{grossman2018extended}
Pinhas Grossman, Scott Morrison, David Penneys, Emily Peters, and Noah Snyder.
\newblock The extended {H}aagerup fusion categories.
\newblock {\em Ann. Sci. \'{E}c. Norm. Sup\'{e}r. (4)}, 56(2):589--664, 2023.

\bibitem[GS16]{grossman2016brauer}
Pinhas Grossman and Noah Snyder.
\newblock The {B}rauer-{P}icard group of the {A}saeda-{H}aagerup fusion
  categories.
\newblock {\em Trans. Amer. Math. Soc.}, 368(4):2289--2331, 2016.

\bibitem[IMP{\etalchar{+}}15]{izumi2015subfactors}
Masaki Izumi, Scott Morrison, David Penneys, Emily Peters, and Noah Snyder.
\newblock Subfactors of index exactly 5.
\newblock {\em Bull. Lond. Math. Soc.}, 47(2):257--269, 2015.

\bibitem[JMS14]{jones2014classification}
Vaughan Jones, Scott Morrison, and Noah Snyder.
\newblock The classification of subfactors of index at most 5.
\newblock {\em Bull. Amer. Math. Soc. (N.S.)}, 51(2):277--327, 2014.

\bibitem[Jon83]{Jones1983}
Vaughan Jones.
\newblock Index for subfactors.
\newblock {\em Invent. Math.}, 72(1):1--25, 1983.

\bibitem[Jon87]{jones1987subfactors}
Vaughan Jones.
\newblock Subfactors of type {${\rm II}_1$} factors and related topics.
\newblock In {\em Proceedings of the {I}nternational {C}ongress of
  {M}athematicians, {V}ol. 1, 2 ({B}erkeley, {C}alif., 1986)}, pages 939--947.
  Amer. Math. Soc., Providence, RI, 1987.

\bibitem[Jon99]{jones1999planar2}
Vaughan F.~R. Jones.
\newblock Planar algebras, i, 1999.

\bibitem[Jon00]{jones2000palgbipartite}
Vaughan Jones.
\newblock The planar algebra of a bipartite graph.
\newblock In {\em Knots in {H}ellas '98 ({D}elphi)}, volume~24 of {\em Ser.
  Knots Everything}, pages 94--117. World Sci. Publ., River Edge, NJ, 2000.

\bibitem[Jon19]{palgnotes}
Vaughan Jones.
\newblock Notes from {P}lanar {A}lgebras course.
\newblock \url{https://my.vanderbilt.edu/jonesvf/files/2020/10/pl21.2019.pdf},
  Fall 2019.

\bibitem[Jon22]{jones1999planar}
Vaughan Jones.
\newblock Planar algebras, {I}.
\newblock {\em New Zealand J. Math.}, 52:1--107, 2021 [2021--2022].

\bibitem[JP11]{jones2011embedding}
Vaughan Jones and David Penneys.
\newblock The embedding theorem for finite depth subfactor planar algebras.
\newblock {\em Quantum Topol.}, 2(3):301--337, 2011.

\bibitem[JS97]{jones1997introduction}
Vaughan Jones and Viakalathur~S. Sunder.
\newblock {\em Introduction to subfactors}, volume 234.
\newblock Cambridge University Press, 1997.

\bibitem[Kaw23]{kawahigashi2023characterization}
Yasuyuki Kawahigashi.
\newblock A characterization of a finite-dimensional commuting square producing
  a subfactor of finite depth.
\newblock {\em Int. Math. Res. Not. IMRN}, (10):8419--8433, 2023.

\bibitem[KS04]{kodiyalam2004jones}
Vijay Kodiyalam and Viakalathur~S. Sunder.
\newblock On {J}ones' planar algebras.
\newblock {\em J. Knot Theory Ramifications}, 13(2):219--247, 2004.

\bibitem[KS09]{kodiyalam2009subfactor}
Vijay Kodiyalam and Viakalathur~S. Sunder.
\newblock From subfactor planar algebras to subfactors.
\newblock {\em Internat. J. Math.}, 20(10):1207--1231, 2009.

\bibitem[Lon94]{longo1994duality}
Roberto Longo.
\newblock A duality for {H}opf algebras and for subfactors. {I}.
\newblock {\em Comm. Math. Phys.}, 159(1):133--150, 1994.

\bibitem[Mon23]{montgomery2023}
Michael Montgomery.
\newblock Colored planar algebras for commuting squares and applications to
  {H}adamard subfactors.
\newblock {\em J. Funct. Anal.}, 285(5):Paper No. 110010, 56, 2023.

\bibitem[MW10]{morrisonwalker2010}
Scott Morrison and Kevin Walker.
\newblock The graph planar algebra embedding theorem, 2010.
\newblock preprint, http://tqft.net/gpa.

\bibitem[Ocn88]{ocneanu1988quantized}
Adrian Ocneanu.
\newblock Quantized groups, string algebras and {G}alois theory for algebras.
\newblock In {\em Operator algebras and applications, {V}ol. 2}, volume 136 of
  {\em London Math. Soc. Lecture Note Ser.}, pages 119--172. Cambridge Univ.
  Press, Cambridge, 1988.

\bibitem[Pet10]{peters2010planar}
Emily Peters.
\newblock A planar algebra construction of the {H}aagerup subfactor.
\newblock {\em Internat. J. Math.}, 21(8):987--1045, 2010.

\bibitem[Pop93]{popa1993markov}
Sorin Popa.
\newblock Markov traces on universal {J}ones algebras and subfactors of finite
  index.
\newblock {\em Invent. Math.}, 111(2):375--405, 1993.

\bibitem[Pop94]{popa1994classification}
Sorin Popa.
\newblock Classification of amenable subfactors of type {II}.
\newblock {\em Acta Math.}, 172(2):163--255, 1994.

\bibitem[Pop95]{popa1995axiomatization}
Sorin Popa.
\newblock An axiomatization of the lattice of higher relative commutants of a
  subfactor.
\newblock {\em Invent. Math.}, 120(3):427--445, 1995.

\bibitem[PP86]{pimsner1986entropy}
Mihai Pimsner and Sorin Popa.
\newblock Entropy and index for subfactors.
\newblock {\em Ann. Sci. \'{E}cole Norm. Sup. (4)}, 19(1):57--106, 1986.

\bibitem[Sch90]{schou2013commuting}
John~K. Schou.
\newblock {\em Commuting squares and index for subfactors}.
\newblock 1990.
\newblock Thesis (Ph.D.)--Odense Universitet Institut for Matematik og
  Datalogi, see also arXiv:1304.5907v1.

\bibitem[Sto21]{stojanovic2021}
Hrvoje Stojanovic.
\newblock {\em New Examples of Irreducible Subfactors of the Hyperfinite {${\rm
  II}_1$} Factor with Rational, Non-Integer Index}.
\newblock 2021.
\newblock Dissertation -- Vanderbilt University, 65 pages.

\bibitem[Wat90]{watatani1990index}
Yasuo Watatani.
\newblock {\em Index for {$C^*$}-subalgebras}, volume~83.
\newblock 1990.

\bibitem[Wen88]{wenzl1988hecke}
Hans Wenzl.
\newblock Hecke algebras of type {$A_n$} and subfactors.
\newblock {\em Invent. Math.}, 92(2):349--383, 1988.

\end{thebibliography}

\end{document}